\documentclass[%
  a4paper,
  onecolumn,
  colorlinks,
]{preprint}

\usepackage[english]{babel}

\usepackage{amsmath}
\usepackage{amssymb}
\usepackage{amsthm}
\usepackage{mathtools}
\usepackage{stmaryrd}
\SetSymbolFont{stmry}{bold}{U}{stmry}{m}{n}
\usepackage{mathrsfs}

\usepackage{graphicx}
\usepackage{xcolor}
\usepackage{subcaption}
\usepackage{tikz}
\usepackage{pgfplots}
\usepackage{pgfplotstable}
\usepackage[linesnumbered, lined, algoruled]{algorithm2e}
\usepackage[T1]{fontenc}
\SetKw{KwBreak}{break}
\usepackage{booktabs}

\usepackage{todonotes}
\makeatletter
\DeclareOldFontCommand{\rm}{\normalfont\rmfamily}{\mathrm}
\makeatother
\newcommand{\ignore}[1]{}

\pgfplotsset{compat = 1.18}
\usepgfplotslibrary{external}
\usetikzlibrary{matrix}

\tikzset{external/system call = {%
    pdflatex \tikzexternalcheckshellescape
    -halt-on-error
    -interaction=batchmode
    -jobname "\image" "\texsource"}}
\tikzexternaldisable

\theoremstyle{plain}

\theoremstyle{definition}

\newtheorem{definition}{Definition}

\newcommand{\R}{\ensuremath{\mathbb{R}}}

\newcommand{\opK}{\ensuremath{\mathscr{K}}}
\newcommand{\opL}{\ensuremath{\mathscr{L}}}
\newcommand{\opM}{\ensuremath{\mathscr{M}}}

\newcommand{\opV}{\ensuremath{\mathscr{V}}}

\newcommand{\bA}{\ensuremath{\boldsymbol{A}}}
\newcommand{\bB}{\ensuremath{\boldsymbol{B}}}
\newcommand{\bC}{\ensuremath{\boldsymbol{C}}}

\newcommand{\bF}{\ensuremath{\boldsymbol{F}}}
\newcommand{\bG}{\ensuremath{\boldsymbol{G}}}
\newcommand{\bH}{\ensuremath{\boldsymbol{H}}}
\newcommand{\bI}{\ensuremath{\boldsymbol{I}}}

\newcommand{\bL}{\ensuremath{\boldsymbol{L}}}
\newcommand{\bM}{\ensuremath{\boldsymbol{M}}}

\newcommand{\bQ}{\ensuremath{\boldsymbol{Q}}}
\newcommand{\bR}{\ensuremath{\boldsymbol{R}}}
\newcommand{\bS}{\ensuremath{\boldsymbol{S}}}
\newcommand{\bT}{\ensuremath{\boldsymbol{T}}}
\newcommand{\bU}{\ensuremath{\boldsymbol{U}}}
\newcommand{\bV}{\ensuremath{\boldsymbol{V}}}

\newcommand{\bX}{\ensuremath{\boldsymbol{X}}}
\newcommand{\bY}{\ensuremath{\boldsymbol{Y}}}

\newcommand{\bOmega}{\ensuremath{\boldsymbol{\Omega}}}
\newcommand{\bPsi}{\ensuremath{\boldsymbol{\Psi}}}
\newcommand{\bDelta}{\ensuremath{\boldsymbol{\Delta}}}
\newcommand{\bLambda}{\ensuremath{\boldsymbol{\Lambda}}}

\newcommand{\ba}{\ensuremath{\boldsymbol{a}}}
\newcommand{\bb}{\ensuremath{\boldsymbol{b}}}
\newcommand{\bc}{\ensuremath{\boldsymbol{c}}}
\newcommand{\be}{\ensuremath{\boldsymbol{e}}}

\newcommand{\br}{\ensuremath{\boldsymbol{r}}}
\newcommand{\bs}{\ensuremath{\boldsymbol{s}}}

\newcommand{\bv}{\ensuremath{\boldsymbol{v}}}

\newcommand{\bx}{\ensuremath{\boldsymbol{x}}}
\newcommand{\by}{\ensuremath{\boldsymbol{y}}}

\newcommand{\bNabla}{\ensuremath{\boldsymbol{\nabla}}}
\newcommand{\bpsi}{\ensuremath{\boldsymbol{\psi}}}

\newcommand{\bTB}{\ensuremath{\boldsymbol{\mathcal{B}}}}
\newcommand{\bTC}{\ensuremath{\boldsymbol{\mathcal{C}}}}
\newcommand{\bTG}{\ensuremath{\boldsymbol{\mathcal{G}}}}

\newcommand{\bTV}{\ensuremath{\boldsymbol{\mathcal{V}}}}
\newcommand{\bTW}{\ensuremath{\boldsymbol{\mathcal{W}}}}
\newcommand{\bTX}{\ensuremath{\boldsymbol{\mathcal{X}}}}
\newcommand{\bTY}{\ensuremath{\boldsymbol{\mathcal{Y}}}}
\newcommand{\bTZ}{\ensuremath{\boldsymbol{\mathcal{Z}}}}

\newcommand{\trans}{\ensuremath{\mkern-1.5mu\mathsf{T}}}
\newcommand{\krp}{\odot}

\newcommand{\fro}{\operatorname{F}}

\DeclareMathOperator{\mdiag}{diag}
\DeclareMathOperator*{\argmin}{argmin}

\newcommand{\plotfontsize}{\small}

\definecolor{matlabblue}{HTML}{0072BD}
\definecolor{matlaborange}{HTML}{D95319}
\definecolor{matlabyellow}{HTML}{EDB120}
\definecolor{matlabpurple}{HTML}{7E2F8E}
\definecolor{matlabgreen}{HTML}{77AC30}
\definecolor{matlablightblue}{HTML}{4DBEEE}
\definecolor{matlabred}{HTML}{A2142F}

\tikzstyle{GMRES} = [
  matlabblue,
  solid,
  line width = 1.5pt,
  mark options = {solid},
  mark         = triangle*
]

\tikzstyle{sGMRES} = [
  matlabgreen,
  solid,
  line width   = 1.5pt,
  mark options = {solid},
  mark         = *
]

\tikzstyle{MLNsGMRES} = [
  matlaborange,
  solid,
  line width   = 1.5pt,
  mark options = {solid},
  mark         = square*
]

\tikzstyle{MLNsGMRESMemoryEfficient} = [
  matlabpurple,
  solid,
  line width   = 1.5pt,
  mark options = {solid},
  mark         = diamond*   
]

\begin{document}


\title{Randomized Tucker-Sketched GMRES}

\author[1]{\mbox{Alberto Bucci}}
\author[2]{\mbox{Martina Iannacito}}
\author[3]{\mbox{Mirjeta Pasha}}
\author[3]{\mbox{Rudi Smith}}

\affil[1]{School of Mathematics, The University of Edinburgh, Edinburgh, EH9 3FD, UK.
\authorcr \email{abucci2@ed.ac.uk}, \orcid{0009-0002-7550-6025}
}

\affil[2]{Department of Mathematics, University of Bologna, Bologna 40126, IT
\authorcr \email{martina.iannacito@unibo.it}, \orcid{0000-0003-3354-2538}
}

\affil[3]{Department of Mathematics, Virginia Tech, Blacksburg, VA 24061, USA.
\authorcr \email{mpasha@vt.edu}, \orcid{0000-0003-4249-2421}
\authorcr \email{smithrgh@vt.edu}, \orcid{0000-0002-6592-2482}
}

\shorttitle{Randomized Tucker Sketched GMRES}
\shortauthor{A. Bucci, M. Iannacito, M. Pasha, R. Smith}
\shortdate{}
\shortinstitute{}


 
\abstract{                                        
We address the problem of solving large-scale tensor-structured linear systems in the Tucker format. In this setting, standard iterative solvers such as GMRES face a fundamental bottleneck: the multilinear ranks of the Krylov basis vectors grow with the iteration count, leading to rapidly increasing tensor operation costs and memory requirements.
To overcome these challenges, we propose two randomized algorithms within the sketched GMRES framework that replace full Arnoldi orthogonalization with short recurrences. The first, RHOSVD-Tucker sGMRES, uses randomized HOSVD with per-iteration rank selection, providing robustness across a wide range of problems. The second method, MLN-Tucker sGMRES, leverages the multilinear Nystr\"om approximation with a fixed rank, enabling streaming computations; the streamability of the approximation further allows, at no additional cost, a memory-efficient reconstruction of the solution from a compact sketched representation of the Krylov basis.
Both methods outperform standard low-rank Tucker solvers in symmetric and non-symmetric settings. Applied to inverse problems, the low-rank Tucker constraint acts as an implicit regularizer; combined with adaptive projected Tikhonov penalization and automatic regularization parameter selection, the methods yield stable reconstructions.
 }

\novelty{}

\maketitle


\section{Introduction}
\label{sec: Introduction}
We consider large-scale tensor equations of the form
\begin{equation}
\label{eq:AXB:teneq}
   \opL(\bTX) = \bTB,
\end{equation}
where $\bTX, \bTB\in \R^{n_1\times n_2 \times \dots \times n_d}$ are $d$-dimensional tensors and $\opL$ is a linear operator expressible as a sum of separable terms. Denoting $\mathrm{vec}(\bTX)$ as the vectorization of $\bTX$,~\Cref{eq:AXB:teneq} takes the equivalent form
\begin{equation}
\label{eq:AXB}
    \left(\sum_{i=1}^N \bA_{1}^{(i)}\otimes \dots \otimes \bA_d^{(i)}\right)
    \mathrm{vec}(\bTX) = \mathrm{vec}(\bTB),
\end{equation}
where $\bA_j^{(i)}\in \R^{n_j\times n_j}$.
Equations of the form~\eqref{eq:AXB:teneq} arise in a wide range of scientific computing applications, including the discretization of high-dimensional partial differential equations~\cite{Kho15}, quantum many-body simulations and quantum chemistry~\cite{KhoK18}, and uncertainty quantification~\cite{BuiGMetal13b,BuiGMetal13a}, where the underlying operators often admit exact or approximate representations as sums of Kronecker products.
{In many of these applications, the underlying operators admit exact or approximate representations as sums of Kronecker products, while the right-hand side $\bTB$ often has a low-rank tensor representation. Moreover, under suitable spectral assumptions on $\opL$, the solution $\bTX$ may inherit a similar low-rank structure} (see, e.g., \cite{Bac23, BalG13}). {Low-rank tensor structure can also be exploited in applications where it is not naturally present but is instead imposed as a modeling or computational constraint,} such as tensor-on-tensor regression~\cite{AhmRB20, Loc18} and imaging inverse problems~\cite{HanNO06, HanJL21}. {In particular, discrete inverse problems provide an important class of problems of the form~\eqref{eq:AXB}, where $\opL$ represents a discretized forward operator, such as the Radon transform in computed tomography~\cite{HanJL21}, convolution in image deblurring~\cite{HanNO06}, or partial Fourier sampling in MRI~\cite{BroS11}.} The data $\bTB$ are typically noisy, i.e.,
$
\bTB = \opL(\bTX_{\mathrm{true}}) + \boldsymbol{\eta}
= \bTB_{\rm true} + \boldsymbol{\eta},
$
with additive noise $\boldsymbol{\eta}$. In these settings, $\opL$ is severely ill-conditioned, with singular values decaying to zero, so direct solution of~\eqref{eq:AXB:teneq} strongly amplifies noise. Classical iterative methods exhibit \emph{semi-convergence}: the error initially decreases as dominant modes are captured, then increases as noise enters the solution. Regularization is therefore required to obtain a stable approximation of $\bTX_{\mathrm{true}}$ \cite{Han10}.

Krylov-based iterative solvers, such as the Generalized Minimal RESidual (GMRES) method~\cite{SaaS86}, have been widely used {for solving large-scale linear tensor equations} especially in the context of discretized PDEs and imaging inverse problems~\cite{Dol13,KreT11,Han10,CalGR99,SemHKetal14,IslMO26,Bat24}. Because of the \textit{curse of dimensionality} where the number of elements of a tensor as well as the storage consumption grow exponentially with the number of the dimensions \cite{GraKT13,KolB09}, Krylov solvers are not directly used on dense tensors defining the tensor equations, but {are typically combined with low-rank tensor formats to reduce storage and computational costs}. However, this introduces a well-known difficulty: {rank growth during the iterations} \cite{BucPR25}. As Krylov vectors are repeatedly formed through linear combinations and orthogonalization steps, their tensor ranks tend to increase, {progressively diminishing the advantages of the low-rank representation}. To control this, rank truncation is commonly used, albeit at the cost of additional approximation error and algorithmic complexity. Several formats have been developed to represent high-dimensional data efficiently, including the Canonical Polyadic decomposition~\cite{Har70,CarC70}, the Tucker format~\cite{Tuc66}, and hierarchical representations such as Tensor Trains~\cite{Ose09} and Hierarchical Tucker~\cite{Gra10}. For very high-order tensors, hierarchical formats are typically preferred to mitigate the curse of dimensionality, while the Tucker format offers a practical compromise between expressiveness and computational efficiency in moderate-order regimes (e.g., $d=3,4$) with large mode sizes. Although the low-rank format allows the adaptation of Krylov solver to large-scale settings, the storage and orthogonalization of the Krylov basis {remain} a major computational bottleneck. In particular, the memory required to store basis vectors grows linearly with the iteration count, while the cost of orthogonalization grows quadratically, quickly becoming prohibitive for problems with millions of unknowns~\cite{PasSGetal23,PasGSetal24,LanPLetal25}. These limitations have motivated a variety of approaches aimed at reducing the cost of Krylov subspace methods, including recycling strategies~\cite{BurG24,PasSK23}, restarted and memory-limited projection schemes~\cite{LinPGetal25,OkuPKetal25, BucGR26}, structured subspace techniques~\cite{LanPLetal25}, and incomplete orthogonalization variants~\cite{SaaS86}. The latter class {provides the basis for} randomized GMRES-type methods, which replace exact orthogonalization with randomized projections based on oblivious subspace embeddings (OSEs) (see~\Cref{def:OSEs}).

This perspective leads to randomized Krylov solvers that construct approximate bases at significantly reduced cost while preserving convergence properties with high probability. Representative examples include the randomized GMRES framework of~\cite{BalG22} and sketched GMRES (sGMRES)~\cite{NakT24} for the matrix cases.
These two existing randomized GMRES variants differ primarily in their orthogonalization strategies. In~\cite{BalG22}, the Krylov basis is fully orthogonalized within the sketched space, closely following
classical GMRES and therefore inheriting the same rank-growth limitations when adapted to
high-dimensional settings. On the contrary, in~\cite{NakT24}, only a part of the sketched Krylov basis is orthogonalized, thanks to short recurrences, a particularly appealing result applicable to low-rank tensor computations~\cite{BucPR25}.
Additionally, randomization can also be exploited beyond the Krylov solver itself. As observed in~\cite{SmiPGetal26}, sketching techniques can be combined with Kronecker and Khatri--Rao structures to perform Tucker tensor operations, such as tensor summation and compression, without forming dense intermediate quantities. While this substantially reduces the cost of individual tensor operations, it does not eliminate the fundamental challenge posed by rank accumulation across the Krylov basis.
Building on these premises, we focus on moderate-order Tucker tensors and investigate randomized strategies for efficiently adapting sketched GMRES (sGMRES)~\cite{NakT24} to this setting.

\subsection{Contributions}
To address the limitations of existing tensor-based solvers, we propose two Tucker-based variants of sGMRES. 
\begin{enumerate}
    \item A first solver, RHOSVD-Tucker sGMRES, builds upon the framework of~\cite{SmiPGetal26} and integrates randomized Higher-Order Singular Value Decomposition (RHOSVD)~\cite{AhmAAetal21} together with an effective heuristic for rank estimation into the Krylov iteration, yielding a more efficient implementation.

    \item The second solver, MLN-Tucker sGMRES, relies on a streamable Tucker compression strategy, the multilinear Nystr\"om (MLN) approximation~\cite{SunGLetal20, BucR24}. At the price of fixing a maximum rank \emph{a priori}, the MLN variant enables efficient operator application, tensor summation, and solution reconstruction. Furthermore, for memory-constrained settings, we present an optional sketching strategy, adapting to the Tucker setting the approach of~\cite{BucPR25}, that replaces the full Krylov basis with a compact sketched representation. This option is broadly applicable whenever the solution is numerically low-rank, but integrates particularly naturally into MLN-Tucker sGMRES, where it incurs no additional cost. 
\end{enumerate}
 
While the fixed-rank constraint imposed by the MLN approach may alter the convergence behavior, it also acts as a form of implicit regularization by restricting the iterates to a manifold of bounded multilinear rank and thereby suppressing high-frequency noise. To more explicitly address ill-posedness, we incorporate an \emph{adaptive projected Tikhonov} regularization strategy~\cite{CalGR99,GazNR15} within each iteration. The regularization parameter is selected automatically via generalized cross-validation (GCV)~\cite{GolHW79}, eliminating the need for oracle-based stopping criteria and enabling stable reconstructions even when the noise level is unknown. To the best of our knowledge, this is the first work demonstrating that compressed tensor decompositions can serve as an effective regularizer in the context of iterative methods for inverse problems. The interplay between explicit Tikhonov penalization and implicit low-rank truncation makes the proposed framework particularly well suited for large-scale, ill-posed tensor inverse problems, while the sGMRES architecture simultaneously addresses the memory and orthogonalization bottleneck described above. 

To demonstrate the effectiveness of our solvers, we consider three benchmark problems. First, a 3D Poisson equation: since our construction recovers short recurrences, the per-iteration cost is greatly reduced, and we show that our approach is competitive with, and in fact faster and more robust than, state-of-the-art CG-based Tucker methods~\cite{IanPS26} which are tailored to symmetric problems and enjoy short recurrences naturally. Second, a 3D convection-diffusion problem, a classical and challenging non-symmetric benchmark. Finally, a 3D image deblurring inverse problem, which illustrates how the low-rank Tucker constraint simultaneously acts as an implicit regularizer and renders the problem computationally tractable.

\paragraph{Outline.}
The remainder of the paper is organized as follows.~\Cref{sec:prelim} reviews tensor notation and the Tucker decomposition, including randomized approximation techniques, GMRES and its low-rank variants such as RHOSVD-Tucker GMRES, which will serve as benchmark for the numerical experiments.~\Cref{sec: Algorithms} presents the two proposed algorithms: RHOSVD-Tucker sGMRES and MLN-Tucker sGMRES. 
~\Cref{sec: Results} reports numerical experiments on a preconditioned 3D Poisson equation, a 3D convection--diffusion problem, and a 3D image deblurring inverse problem. Finally,~\Cref{sec:conclusions} provides concluding remarks.

\section{Preliminaries}
\label{sec:prelim}

\subsection{Tensor Notation and Algebraic Operations}
\label{sec:Preliminaries}

We adopt standard notation whereby lowercase italic letters ($x$) represent scalars, bold lowercase letters ($\bx$) denote vectors, bold uppercase letters ($\bX$) signify matrices, bold calligraphic letters ($\bTX$) identify tensors, and script letters are used for operators ($\opL$) or sets of matrices ($\mathscr{X},\mathscr{Y}$). An order $d$ tensor is defined as $\bTX \in \R^{n_{1} \times n_{2} \times \dots \times n_{d}}$, where individual elements are accessed through $d$ indices $x_{i_{1} i_{2} \dots i_{d}}$ with $1 \le i_{k} \le n_{k}$.

\subsection{Tensor Operations}
\label{sec:TensorOperations}

We rely on several foundational operations from tensor algebra~\cite{BalK25}. The construction and manipulation of low-rank tensors depends on structured matrix products, which allow us to relate high-dimensional data to lower-dimensional factor matrices.
\begin{definition}[Khatri-Rao Products]

For two matrices $\bA \in \R^{I \times r}$ and $\bB \in \R^{J \times r}$, the Khatri-Rao product, denoted by $\bA \krp \bB \in \R^{IJ \times r}$, is defined as the column-wise Kronecker product
$$\bA \krp \bB = \begin{bmatrix} \ba_{1} \otimes \bb_{1} & \ba_{2} \otimes \bb_{2} & \dots & \ba_{r} \otimes \bb_{r} \end{bmatrix}$$
where $\otimes$ denotes the Kronecker product.
\end{definition}
We also need to project or transform a tensor along a specific dimensions using the tensor-times-matrix and matricized tensor times Khatri-Rao products.
\begin{definition}[Tensor-Times-Matrix (\texttt{TTM})]
\label{def:ttm}
The mode-$k$ tensor-times-matrix product (\texttt{TTM}) of a tensor $\bTX \in \R^{n_1 \times \dots \times n_d}$ with a matrix $\bU \in \R^{J \times n_k}$ is denoted by $\bTY = \bTX \times_k \bU$. This operation applies the linear transformation $\bU$ to each mode-$k$ fiber of $\bTX$, resulting in a new tensor $\bTY \in \R^{n_1 \times \dots \times n_{k-1} \times J \times n_{k+1} \times \dots \times n_d}$. Computationally, the \texttt{TTM} is used in its matricized form, where the mode-$k$ product corresponds to standard matrix multiplication of some unfolding. Denoting the mode-$k$ unfolding of a tensor by the subscript $(k)$, this is expressed as
$$\bY_{(k)} = \bU \bX_{(k)}.$$
To denote a sequence of successive \texttt{TTM}s across multiple modes, we use the product
\begin{equation}
    \bTX \times_1 \bU_1 \times_2 \bU_2 \times_3 \bU_3 := \bTX \times_{k=1}^{3}\bU_{k}.
\end{equation}
\end{definition}
To contract a tensor with multiple factor matrices under the Khatri-Rao umbrella, we use the matricized tensor times Khatri-Rao product.
\begin{definition}[Matricized Tensor Times Khatri-Rao Product (\texttt{MTTKRP})]
\label{def:mttkrp}
Given an order-$d$ tensor $\bTX \in \R^{n_1 \times \dots \times n_d}$ and a sequence of factor matrices $\{\bU_i \in \R^{n_i \times r}\}_{i=1}^d$, the \texttt{MTTKRP} with respect to mode $k$ produces a matrix $\bM^{(k)} \in \R^{n_k \times r}$. It is formulated as the matrix multiplication of the mode-$k$ unfolding of $\bTX$ and the Khatri-Rao product of all factor matrices except the $k$-th one 
\begin{equation}
\label{eq:mttkrp}
    \bM^{(k)} = \bX_{(k)} \left( \bU_d \krp \dots \krp \bU_{k+1} \krp \bU_{k-1} \krp \dots \krp \bU_1 \right).
\end{equation}
and is denoted $\texttt{MTTKRP}(\bTX; \{\bU_j\}_j, k)$.
\end{definition}

\subsection{The Tucker Decomposition}
\label{sec:TuckerTensors}
The Tucker decomposition~\cite{Tuc66} is a low rank decomposition formed by projecting an order-$d$ tensor onto smaller, mode-specific latent spaces through the \texttt{TTM} operation in~\Cref{def:ttm}, yielding a $d$ dimensional tensor core and a set of $d$ factor matrices. Algorithms have been made to operate on the compact components of the Tucker format, however, the limitations of such algebraic operations usually lie with respect to the core. 

More formally, using~\Cref{def:ttm}, an order-$d$ tensor $\bTX\in\R^{n_1\times\cdots\times n_d}$ with a specified multilinear rank $(r_1,\dots,r_d)$ is factorized into a core tensor $\bTG\in\R^{r_1\times\cdots\times r_d}$ and a collection of factor matrices $\{\bU_k\in\R^{n_k\times r_k}\}_{k=1}^d$. The Tucker approximation is constructed through a multi \texttt{TTM} through all the $k$ modes 
\begin{equation}
\label{eq:TuckerFormat}
    \bTX \approx \bTG \times_1 \bU_1 \times_2 \bU_2\cdots\times_d \bU_d,
\end{equation}
although we emphasize that this is rarely formed, but rather stored as a core and factors.
While the factor matrices compress each mode from a dimension of $n_k$ to $r_k$, the Tucker format still has its weakness through the core scaling requiring $\mathcal{O}(r^d)$ units of memory. Consequently while valuable for low-order data, the Tucker format becomes intractable quickly, so usually sees its applications for modest dimension sizes. Throughout this paper, we operate under the assumption that the compressed representation of $\bTB$ from~\Cref{eq:AXB:teneq} has been given \textit{a priori} in Tucker format, and that the order $d$ is sufficiently small to permit algorithmic manipulation of the core $\bTG$. 

Several deterministic and randomized algorithms for computing Tucker approximations have been proposed in the literature~\cite{KolB09,PeaM25, CheW19,CheWY21, IanPPetal26}. In this work we primarily employ the randomized HOSVD~\cite{AhmAAetal21} and the Multilinear Nystr\"om method~\cite{SunGLetal20,BucR24}, although other Tucker compression schemes, such as sequential and parallel variants~\cite{MinSK20, BucH25, HasN25, MinLB24} could be used as well. These methods are employed to handle the operator applications and linear combinations arising in the GMRES orthogonalization phase and the reconstruction of the final solution.

\subsection{GMRES and variants for linear systems}\label{sec:prelim:GMRES}

The GMRES iterative solver~\cite{SaaS86} constructs an approximate solution to the linear system $\bA\bx = \bb$ within a Krylov subspace of growing dimension. In particular, starting from an initial guess $\bx_0$, the approximate solution at the $k$-th iteration is given by
\begin{equation}
    \bx_k = \bx_0 + \bV_k\by_k,
\end{equation}
where $\bV_k\in\R^{n\times k}$ forms a basis for the Krylov subspace $\opK_k(\bA, \br_0) = {\rm span}\bigl\{\br_0, \bA\br_0, \dots, \bA^{k-1}\br_0\bigr\}$, with $\br_0 = \bb - \bA\bx_0$ denoting the initial residual. The vector $\by_k \in\R^{k}$ is determined by solving the least-squares problem
\begin{equation}
\label{eq:LSGMRES:y}
    \by_k =\arg\min_{\by\in\R^{k}} \|\bA\bV_k\by - \br_0\|_2.
\end{equation}
In practice, if the basis for $\opK_k$ is constructed using a full orthogonalization of its elements via the Arnoldi process, the resulting matrix $\bV_k$ has orthonormal columns. We then rely on the Arnoldi relation
\begin{equation}
\label{eq:Arnoldi}
    \bA\bV_k = \bV_{k+1}\overline{\bH}_k,
\end{equation}
where $\overline{\bH}_k\in\R^{(k+1)\times k}$ is an upper-Hessenberg matrix containing the orthogonalization coefficients. We can then rewrite~\eqref{eq:LSGMRES:y} in a simplified manner as
\begin{equation}
\label{eq:LS:y=Hy-betae1}
    \by_k = \arg\min_{\by\in\R^{k}}\|\overline{\bH}_k\by - \beta\be_1\|_2,
\end{equation}
where $\beta = \|\br_0\|_2$ and $\be_1 = (1,0,\dots,0)^{\trans}\in\R^{k+1}$. 
The convergence of GMRES is sensitive to the orthogonality of the Krylov basis $\bV_k$~\cite{PaiRS06}. However, maintaining strict orthogonality across all basis vectors introduces a computational and memory bottleneck. To address this, various modifications to the classical GMRES algorithm have been proposed. Among the most widely adopted strategies to accelerate the solver and reduce costs are restarting~\cite{Mor02, Saa03b}, preconditioning~\cite{Saa03b}, and incomplete orthogonalization~\cite{SaaS86}, the latter of which serves as the foundation for randomized GMRES variants. In the incomplete orthogonalization approach, the newly generated basis vector $\bv_k$ is orthogonalized only against the $k_{\text{trunc}}$ most recently computed basis vectors, reducing the computational overhead of the orthogonalization phase. Because of this, while the Arnoldi relation~\eqref{eq:Arnoldi} remains satisfied and $\by_k$ can still be computed via~\eqref{eq:LS:y=Hy-betae1}, the basis matrix $\bV_k$ loses its strict orthogonality. Consequently, the resulting $\by_k$ no longer represents the exact mathematical solution to~\eqref{eq:LSGMRES:y}. As such, the cost per iteration is lowered but at the expense of delayed convergence. 
Motivated by reducing orthogonalization costs, randomized variants of GMRES have also been proposed~\cite{Woo14,MarT20,MurDMetal23} which rely on the use of oblivious subspace embeddings (OSEs). By implicitly embedding the high-dimensional Arnoldi process into a low-dimensional sketched space via an OSE, one can construct the Krylov basis and solve the corresponding least-squares problem~\eqref{eq:LSGMRES:y} at a fraction of the cost without sacrificing accuracy.

\begin{definition}[Oblivious Subspace Embedding]\label{def:OSEs}
Let $\varepsilon \in (0,1)$ and $\delta \in (0,1)$. A random matrix $\bS\in \mathbb{R}^{s\times n}$ is said to satisfy the $(\varepsilon, \delta, k)$-OSE property if, for any fixed $k$-dimensional subspace $\mathcal{V} \subseteq \mathbb{R}^n$, the following holds with probability at least $1-\delta$
\begin{equation}
    (1-\varepsilon)\|\bx\|_2^2 \leq \|\bS \bx\|_2^2 \leq (1+\varepsilon)\|\bx\|_2^2 \quad \forall \bx \in \mathcal{V}.
\end{equation}
\end{definition}

Since the set of all possible residual vectors forms a low-dimensional subspace, the OSE guarantees that minimizing the sketched norm $\|\bS \bx\|_2$ yields a near-optimal approximation to the true high-dimensional minimum. Exploiting this to reduce orthogonalization costs,~\cite{NakT24} introduced the sketched GMRES (sGMRES) method. By projecting the high-dimensional Arnoldi process through an OSE, sGMRES avoids solving the exact least-squares problem~\eqref{eq:LSGMRES:y}. Instead, it derives an approximate solution by solving the reduced problem
\begin{equation}
\label{eq:LSsGMRES:y}
    \by_k = \arg\min_{\by\in\R^{k}}\|\bS \bA\bV_k\by - \bS\br_0\|_2.
\end{equation}

\subsubsection{Low rank variants}
\label{sec:lowrankvariants}

Several low-rank tensor variants of GMRES have been proposed~\cite{KreT11, BalG13, Dol13, MenAC25}, both for classical linear equations and tensor multilinear ones. At iteration $k$, the tensor-adapted GMRES method constructs a tensor $\bTV_k$ to expand a basis of the Krylov tensor subspace
$\opK_k(\opL, \bTB) \coloneqq \mathrm{span}\{\bTB, \opL(\bTB), \dots, \opL^{k-1}(\bTB)\}$. This construction is carried out in a chosen low-rank tensor format.
The orthonormalization coefficients, computed via tensor inner products, are used to assemble the upper-Hessenberg matrix $\overline{\bH}_k$, which in turn defines the update vector $\by_k$ as in~\eqref{eq:LSGMRES:y}.

In all these low-rank GMRES adaptations, truncation at prescribed rank or accuracy is employed to prevent memory errors. 
Indeed, the computation of the new Krylov basis element during the Arnoldi iteration, 
and the orthogonalization against previous basis elements causes the rank to grow as in the matrix case, but with more severe memory usage. Without truncation, one would witness exponential rank growth and thus rapid memory usage leading to computational intractability. 

Typically, the SVD, HOSVD, ST-HOSVD~\cite{VanVM12,DonYQetal23,DeLDV00b}, TT-SVD~\cite{Ose11}, or their randomized variants~\cite{AldHetal23, MinSK20, BucH25, BucR24, BucV26, KreVV23}, are used depending on the chosen format. However, because of truncation, the new Krylov basis element does not, in general, lie exactly in the Krylov subspace, so that the resulting basis is only approximate~\cite{SimS03, Dol13}. 

Combinations of GMRES with the Tucker format have been considered in~\cite{KreT11, MenAC25, SmiPGetal26}. We note in the deterministic framework, the summation of Tucker tensors as part of GMRES requires the densification of increasingly large core tensors in each iteration as part of the truncation step mentioned above. To avoid this the Tucker GMRES using KRPSum-Tucker~\cite{SmiPGetal26}, which we herein call RHOSVD-Tucker GMRES was proposed. In particular all necessary Tucker tensor summations are done in a randomized manner via Khatri-Rao or Kronecker sketching using an effective rank estimation. 
While the improved summation handles the arithmetic efficiency of RHOSVD-Tucker GMRES compared to standard deterministic Tucker GMRES, it retains the full orthogonalization requirement of classical GMRES so can be expensive in long iteration schemes. Thus we propose sGMRES derived variants.


\section{The Algorithms}
\label{sec: Algorithms}

In this section, we introduce our novel sketched GMRES algorithms for solving low-rank tensor equations in the Tucker format, RHOSVD-Tucker sGMRES and MLN-Tucker sGMRES, alongside the existing method for baseline comparison RHOSVD-Tucker GMRES.
 
Our first proposed algorithm, RHOSVD-Tucker sGMRES, naturally extends~\cite{SmiPGetal26} by embedding the Khatri-Rao summation framework within the sGMRES architecture and replaces full orthogonalization with a short-recurrence.

Our second proposed algorithm within the sGMRES framework, MLN-Tucker sGMRES, does not rely on RHOSVD techniques but instead approximations via the Multilinear Nystr\"om method. By relying on a streamable Tucker approximation, it avoids storing the full Krylov basis tensors, as new iterates are generated while retaining only their sketches, reducing memory usage and potentially enabling the use of larger ranks or the avoidance of restarting strategies.

\subsection{RHOSVD-Tucker GMRES}
\label{sec: RHOSVD-Tucker GMRES}

To overcome the summation memory issue in the Tucker format,~\cite{SmiPGetal26} uses Khatri-Rao product based randomized sketching, which implicitly rounds the final sum without ever realizing the massive intermediate structure. This has successfully been applied to other tensor formats too such as TT~\cite{AldBGetal25}. Specifically, to approximate a sum of $N$ Tucker tensors $\bTC = \sum_{i=1}^N \bTX^{(i)}$ with cores $\bTG^{(i)}\in\R^{r_1^{(i)}\times\dots\times r_d^{(i)}}$ and factors $\bU_k^{(i)} \in \R^{n_k \times r_k^{(i)}}$, the authors use the \texttt{MTTKRP} (\Cref{def:mttkrp}) to find the sketched mode-$k$ unfolding as
\begin{equation}
    \bY_k = \bC_{(k)} \bS^{(k)} = \sum_{i=1}^N \bU_k^{(i)} \bG_{(k)}^{(i)} \left( \bigodot_{j\neq k} \bM_{d-j+1}^{(i)} \right),
\end{equation}
where $\bM_j^{(i)} = (\bU_j^{(i)})^{\trans} \bS_j \in \R^{r_j^{(i)} \times s}$ is constructed using parts of a global sketching operator $\bS^{(k)} = \bigodot_{j \neq k} \bS_j$ with independent Gaussian matrices $\bS_j \in \R^{n_j \times s}$. Subsequently a low-dimensional orthonormal basis $\widehat{\bU}_k$ for each mode via an economy QR factorization is found and used in a core projection $\widehat{\bTC} = \sum_{i=1}^N \bTG^{(i)} \times_1 (\widehat{\bU}_1^{\trans}\bU_1^{(i)}) \dots \times_d (\widehat{\bU}_d^{\trans}\bU_d^{(i)})$ before finally applying the ST-HOSVD to the intermediate core $\widehat{\bTC}$, yielding a compact Tucker decomposition that satisfies some user-defined tolerance $\epsilon$. 
Instead of using a static sketch size, the authors propose to dynamically determine the uniform sketch dimension $s$ via estimating the effective rank of the sum by analyzing the spectral decay of the energy-weighted concatenated factor matrices $\bV_{k} = \left[ \|\bTG^{(1)}\|_{\fro} \bU_k^{(1)}, \dots, \|\bTG^{(N)}\|_{\fro} \bU_k^{(N)} \right],$ whereby computing the eigenvalues of the Gram matrix $\bV_k^{\trans} \bV_k$, a minimal rank $\widetilde{r}_k$ required to capture its dominant spectral energy up to a relative threshold can be found. The uniform target sketch size is then conservatively set to $s = \max_k (\widetilde{r}_k) + p$, where $p$ is an oversampling parameter. RHOSVD-Tucker GMRES, defined in~\Cref{alg:tucker_gmres}, takes GMRES with Tucker tensor input but allows for the summation steps in the operator application, orthogonalization, and reconstruction of the final solution to be done efficiently using randomization. Herein, we denote the randomized summation technique from~\cite{SmiPGetal26} as \texttt{RoundSum}.

\begin{algorithm}[t]
    \small
    \caption{\texttt{RHOSVD-Tucker GMRES}}
    \label{alg:tucker_gmres}
    \SetKwInOut{Input}{Input}
    \SetKwInOut{Output}{Output}
    \DontPrintSemicolon

    \Input{Linear operator $\opL$; Tucker tensor $\bTB$; tolerance \texttt{tol}; max iterations $N_{\max}$}
    \Output{Approximate solution in Tucker format $\widetilde{\bTX}$, such that $\|\opL(\widetilde{\bTX})- \bTB\|_{\fro}\leq \texttt{tol}$}

    $\beta \gets \|\bTB\|_{\fro}$\;
    $\bTV_1 \gets \bTB / \beta$\;
    Initialize Hessenberg matrix $\bH \in \R^{(N_{\max}+1) \times N_{\max}}$\;
    
    \For{$k \gets 1$ \KwTo $N_{\max}$}{
        $\bTW \gets \bTV_k$\;
        $\bTW \gets \texttt{RoundSum}\left(\opL(\bTW)\right)$\;

        \For{$j \gets 1$ \KwTo $k$}{
            $h_{j,k} \gets \langle \bTV_j, \bTW \rangle_{\fro}$\;
        }
        
        $\bTW \gets \texttt{RoundSum}\left(\bTW - \sum_{j=1}^k h_{j,k}\bTV_j\right)$\; \label{lin:RHOSVD_orthogonalization}
        $h_{k+1, k} \gets \|\bTW\|_{\fro}$\;
        
        \If{$h_{k+1, k} \approx 0$}{
            \KwBreak\;
        }
        
        $\bTV_{k+1} \gets \bTW / h_{k+1, k}$\;
        
        $\by_k \gets \arg\min_{\by} \| \beta \be_1 - \bH_{1:k+1, 1:k} \by \|_2 (+\lambda\|\by\|_2^2$ if regularization is needed)\;
        $\mathrm{res}_k \gets \dfrac{\| \beta \be_1 - \bH_{1:k+1, 1:k} \by_k \|_2}{\beta}$ \label{lin:RHOSVD_residual} \tcp*{Relative residual}
        
        \If{$\mathrm{res}_k < \textup{\texttt{tol}}$}{
            \KwBreak\; \label{lin:residual_check}
        }
    }
    
    $\widetilde{\bTX} \gets \texttt{RoundSum}\left(\sum_{i=1}^k y_i \bTV_i\right)$\;    \label{lin:RHOSVD_final_sum}
    
    
    \Return $\widetilde{\bTX}$\;
\end{algorithm}

\subsection{RHOSVD-Tucker sGMRES}
\label{sec: RHOSVD-Tucker sGMRES}

While~\Cref{alg:tucker_gmres} accelerates individual tensor operations, the cost of full Arnoldi orthogonalization grows with the iteration count and becomes the dominant computational bottleneck which we overcome in sGMRES. More specifically, the primary algorithmic deviations from the baseline RHOSVD-Tucker GMRES are within the orthogonalization and residual tracking phases, specifically at Lines~\ref{lin:RHOSVD_orthogonalization},~\ref{lin:RHOSVD_residual}, and~\ref{lin:residual_check} of~\Cref{alg:tucker_gmres}. Rather than enforcing full orthogonalization, the newly generated basis tensor is orthogonalized strictly against the $k_{\text{trunc}} > 0$ most recently computed basis tensors (see Line~\ref{lin:sRHOSVD_orthogonalization} of~\Cref{alg:tucker_sgmres}).

As previously mentioned, the localized orthogonalization sacrifices orthogonality of the Krylov basis, so~\eqref{eq:LS:y=Hy-betae1} is no longer valid. Instead, the sketched least-squares problem~\eqref{eq:LSsGMRES:y} is targeted (Line~\ref{lin:sketched_residual} of~\Cref{alg:tucker_sgmres}). For algorithmic efficiency, a Khatri-Rao structured Gaussian matrix $\bS$ is drawn at initialization and used to construct the sketched basis throughout the iteration. We also note the sketched residual may underestimate the true residual so a stopping criterion based directly on the prescribed tolerance $\texttt{tol}$ may lead to premature termination. Thus we introduce a correction factor $\eta\in(0,1)$ and terminate the iterations only when the sketched residual falls below $\eta \cdot \texttt{tol}$ (see Line~\ref{lin:sketched_residual_check} in~\Cref{alg:tucker_sgmres}).

\begin{algorithm}[t]
    \small
    \caption{\texttt{RHOSVD-Tucker sGMRES}}
    \label{alg:tucker_sgmres}
    \SetKwInOut{Input}{Input}
    \SetKwInOut{Output}{Output}
    \DontPrintSemicolon

    \Input{Linear operator $\opL$; Tucker tensor $\bTB$; sketch size $s$; tolerance \texttt{tol}; sketched correction factor $\eta$; max iterations $N_{\max}$; truncation window $k_{\text{trunc}}$}
    \Output{Approximate solution in Tucker format $\widetilde{\bTX}$, such that $\|\opL(\widetilde{\bTX})- \bTB\|_{\fro}\leq \texttt{tol}$}

    $\beta \gets \|\bTB\|_{\fro}$\;
    $\bTV_1 \gets \bTB / \beta$\;
    Generate Gaussian matrices $\{\bS_k\}_{k=1}^{d}$, with $\bS_k\in\mathbb{R}^{n_k\times s}$\;
    $\bS \gets \odot_{k=1}^d \bS_k$\;
    $\bb_s \gets \bS^{\trans}\mathrm{vec}(\bTB)$\;\label{lin:sb}
    $\beta_s \gets \|\bb_s\|_2$\;
    $\bM \gets [\,\,]$\tcp*{initialization for   $\bS\bA\bV_k$}
    \For{$k \gets 1$ \KwTo $N_{\max}$}{
        $\bTW \gets \bTV_k$\;
        $\bTW \gets \texttt{RoundSum}\left(\opL(\bTW)\right)$\;\label{lin:roundsumL}
        $\bM \gets [\bM, \bS^{\trans}\mathrm{vec}(\bTW)]$\label{lin:sw}\;

        $\bTW \gets \texttt{RoundSum}\left(\bTW - \sum_{j=\max(1, k-k_{\text{trunc}}+1)}^k \langle \bTV_j, \bTW \rangle_{\fro} \bTV_j\right)$\; \label{lin:sRHOSVD_orthogonalization}
        
        \If{$\|\bTW\|_{\fro} \approx 0$}{
            \KwBreak\;
        }
        
        $\bTV_{k+1} \gets \bTW / \|\bTW\|_{\fro}$\;
        
        $\by_k \gets \arg\min_{\by}\, \left\| \bM \by - \bb_s \right\|_2 (+\lambda\|\by\|_2^2$ if regularization is needed)\;
        $\mathrm{res}_k \gets \dfrac{\|\bM \by_k - \bb_s\|_2}{\beta_s}$ \label{lin:sketched_residual} \tcp*{Sketched relative residual}

        \If{$\mathrm{res}_k < \eta \label{lin:sketched_residual_check} \cdot \textup{\texttt{tol}}$}{
            \KwBreak\;  
        }
    }    
    $\widetilde{\bTX} \gets \texttt{RoundSum}\left(\sum_{i=1}^k y_i \bTV_i\right)$\;\label{lin:roundsum:x}
    
    
    \Return $\widetilde{\bTX}$\;
\end{algorithm}

\subsection{MLN-Tucker sGMRES}
\label{sec: MLN-Tucker sGMRES}

An alternate procedure for computing low-rank approximations in the Tucker format is the MLN algorithm, sometimes favored for its single-pass, streamable properties. Given a tensor $\bTX\in \mathbb{R}^{n_1 \times \dots \times n_d}$, finding a uniform rank $r$ approximation begins by drawing two sets of dimension reduction matrices (DRMs), $\{\bOmega_k\}_{k=1}^d$ and $\{\bPsi_k\}_{k=1}^d$, of sizes $\prod_{j\neq k} n_j \times r$ and $n_k\times (r + p)$, respectively, where $p$ is an oversampling parameter. The algorithm then computes the Tucker approximation via
\begin{equation}
    (\bTX\times_{k=1}^d \bPsi_k^{\trans}) \times_{k=1}^d \bX_{(k)}\bOmega_k(\bPsi_k^{\trans}\bX_{(k)}\bOmega_k)^\dagger,
\end{equation}
where $\dagger$ denotes the Moore-Penrose pseudoinverse.
When the tensor $\bTX$ is already given in the Tucker format~\eqref{eq:TuckerFormat} the algorithm can be significantly accelerated by employing Khatri-Rao structured sketching matrices in place of $\bOmega_k$. We refer the reader to~\cite{CamEMetal25} for a detailed analysis of the associated theoretical properties.
As for the oversampling parameter $p$, one uses it for both improving the accuracy and numerical stability of the method~\cite{BucR24}. Typically, independent random matrices are drawn to perform the left and right sketching phases, however to recycle computations and reduce the overall cost, our algorithm instead uses a subset of the left sketching matrices to construct the right ones. 
Specifically, we let $\overline{\bPsi}_k$ denote the matrix consisting of the first $r$ columns of $\bPsi_k$, and subsequently set $\bOmega_k$ as the Khatri-Rao product of these submatrices
\begin{equation}
    \bOmega_k = \overline{\bPsi}_d \odot \dots \odot \overline{\bPsi}_{k+1} \odot \overline{\bPsi}_{k-1} \odot \dots \odot \overline{\bPsi}_1.
\end{equation}
We now outline the process of computing a low-rank approximation of a linear combination of Tucker tensors using the MLN framework equipped with Khatri-Rao structured sketching operator.
Importantly in this algorithm each constituent tensor $\bTX^{(j)}$ is accessed exactly once, and only its compact sketched representation is accumulated in memory. Formally, let $\{\bTX^{(j)}\}_{j=1}^N$ be a collection of $N$ tensors of order $d$, where each $\bTX^{(j)} \in \R^{n_1 \times \dots \times n_d}$ is represented by its Tucker decomposition $\bTX^{(j)} = \bTG^{(j)} \times_{k=1}^d \bU_k^{(j)}$.
Given a scalar coefficient vector $\bc = [c_1, \dots, c_N]^{\trans} \in \R^{N}$, our objective is to construct an efficient MLN approximation for the exact linear combination $\widetilde{\bTX} = \widetilde{\bTG}\times_{k=1}^{d}\widetilde{\bU}_k \approx \sum_{j=1}^N c_j \bTX^{(j)}.$ For each summand $\bTX^{(j)}$ we first project the factor matrices onto these sketching operators
\begin{equation}
    \bU_{\bPsi k}^{(j)} := \bPsi_k^{\trans} \bU_k^{(j)} \in \R^{(r+p) \times r^{(j)}}, \qquad \overline{\bU}_{\bPsi k}^{(j)} := \overline{\bPsi}_k^{\trans} \bU_k^{(j)} \in \R^{r \times r_k^{(j)}}.
\end{equation}
Because $\overline{\bPsi}_k$ is strictly a block partition of $\bPsi_k$, the matrix $\overline{\bU}_{\bPsi k}^{(j)}$ does not need to be explicitly computed from scratch; it is extracted as the first $r$ rows of $\bU_{\bPsi k}^{(j)}$, thereby bypassing any unnecessary dense matrix multiplications. Next we construct the corresponding MLN sketched cores and factor contributions for each summand. Because the linear projection distributes over the Tucker format, the $j$-th sketched core tensor $\bTG_{\Psi}^{(j)}$ is computed via small tensor-times-matrix products $\bTG_{\Psi}^{(j)} = \bTX^{(j)} \times_{k=1}^d \bPsi_k^{\trans} = \bTG^{(j)} \times_{k=1}^d \bU_{\bPsi k}^{(j)}.$ Simultaneously, the mode-$k$ factor contributions $\bF_k^{(j)}$ are obtained by applying the \texttt{MTTKRP} from~\Cref{def:mttkrp}. Applying the \texttt{MTTKRP} to the small core tensor and the truncated sketched factors yields
\begin{equation}
    \bF_{\bPsi k}^{(j)} = \bU_k^{(j)} \, \texttt{MTTKRP}\!\left(\bTG^{(j)}, \{\overline{\bU}_{\bPsi \ell}^{(j)}\}_{\ell}, k\right).
\end{equation} 
By isolating the \texttt{MTTKRP} operation to the core tensor, we avoid forming the Khatri-Rao product of the full-dimension factor matrices. Once these summand-specific quantities have been computed, we sum them to form the global linear combinations. Given the scalar coefficients $c_j$, the intermediate sketched core $\bTG_{\bPsi}$ and the accumulated factor contributions $\bF_{\bPsi k}$ are assembled
\begin{equation}
\label{eqn: recoveryinputs}
    \bTG_{\bPsi} = \sum_{j=1}^N c_j \bTG_{\bPsi}^{(j)}, \qquad \bF_{\bPsi k} = \sum_{j=1}^N c_j \bF_{\bPsi k}^{(j)}.
\end{equation}
The exact MLN approximation is now equivalent to the Tucker tensor $\bTG_{\bPsi} \times _{k=1}^d \bF_{\bPsi k}(\bPsi_k^{\trans}\bF_{\bPsi k})^\dagger$. However, because $\bTG_{\bPsi}$ has uniform sizes $r+p$, this representation suffers from unnecessarily inflated multilinear ranks, taking on the dimension of the sketch sizes $r+p$ rather than the desired target ranks $r$. To reduce the representation back to rank $r$, we use a QR compression step. Specifically, we compute the factorization
\begin{equation}
    \bQ_k \bR_k = \bPsi_k^{\trans}\bF_{\bPsi k},
\end{equation}
where $\bQ_k \in \R^{(r+p) \times r}$ has orthonormal columns and $\bR_k \in \R^{r \times r}$ is upper-triangular. Substituting this factorization into the pseudoinverse yields $(\bPsi_k^{\trans}\bF_k)^\dagger = \bR_k^{\dagger} \bQ_k^{\trans}$. We can then absorb the orthogonal matrix $\bQ_k^{\trans}$ directly into the core tensor while applying the pseudoinverse of the upper-triangular factor $\bR_k$ to the factor matrices. The final, rank-compressed core tensor $\widetilde{\bTG}$ and factor matrices $\widetilde{\bU}_k$ are thus computed as
\begin{equation}
    \widetilde{\bTG} = \bTG_{\bPsi} \times_{k=1}^d \bQ_k^{\trans}, \qquad \widetilde{\bU}_k = \bF_{\bPsi k} \bR_k^{\dagger}.
\end{equation}
The overall procedure can be naturally divided into three phases: a sketching phase (\texttt{MLN\_Sketch},~\Cref{alg:mln_sketch}), a summation phase (\texttt{MLN\_Sum},~\Cref{alg:mln_sum}), in which the sketched contributions are accumulated in a streaming fashion, and a recovery phase (\texttt{MLN\_Recovery},~\Cref{alg:mln_recovery}). We further define \texttt{MLN\_RoundSum} as the sequential application of these three routines. 
The full pseudocode for MLN-Tucker GMRES is outlined in~\Cref{alg:mln_tucker_sgmres}. Of particular note is a novelty in Lines~\ref{lin: mln_sketchB} and~\ref{lin: mln_diag} of~\Cref{alg:mln_tucker_sgmres} where normally we would have to generate an additional sketching matrix and compute an additional sketch as in RHOSVD-Tucker sGMRES, Lines~\ref{lin:sb} and~\ref{lin:sw} of~\Cref{alg:tucker_sgmres}. Instead from \texttt{MLN\_Sketch}~\Cref{alg:mln_sketch} we can receive an optional output,  $\mdiag(\bTG_{\bPsi})$, the superdiagonal of $\bTG_{\bPsi}$, which is implicitly computed and reuse this as the sketch in these lines. To see why, let the exact linear combination be $\bTX = \sum_{j=1}^N c_j \bTX^{(j)}$, such that the global sketched core is defined as $\bTG_{\bPsi} = \bTX \times_{k=1}^d \bPsi_k^{\trans}$. The $i$-th scalar element of its superdiagonal, $[\mdiag(\bTG_{\bPsi})]_i$, corresponds to the simultaneous projection of $\bTX$ exclusively onto the $i$-th row of each transposed sketching matrix. Letting $\bpsi_{k,i}$ denote the $i$-th column of $\bPsi_k$, this element evaluates to 
\begin{equation}
     [\mdiag(\bTG_{\bPsi})]_i = (\bTX \times_{k=1}^d \bpsi_{k,i}^{\trans}).
\end{equation}
Applying the standard identity connecting multilinear products to vectorization and the Kronecker product, this scalar operation can be flattened into an exact inner product
\begin{equation}
    [\mdiag(\bTG_{\bPsi})]_i = (\bpsi_{1,i} \otimes \dots \otimes \bpsi_{d,i})^{\trans} \mathrm{vec}(\bTX).
\end{equation}
 By definition, the Kronecker product of the individual $i$-th columns, $\bpsi_{1,i} \otimes \dots \otimes \bpsi_{d,i}$, forms precisely the $i$-th column of the full Khatri-Rao product $\bPsi_1 \odot \dots \odot \bPsi_d$. Stacking these scalar components for all $i$ along the superdiagonal reconstructs the full matrix-vector operation, yielding the final equivalence $\mdiag(\bTG_{\bPsi}) = (\bPsi_1 \odot \dots \odot \bPsi_d)^{\trans} \mathrm{vec}\left(\sum_{j=1}^N c_j\bTX^{(j)}\right).$
The MLN-Tucker sGMRES algorithm follows the same overall structure as RHOSVD-Tucker sGMRES, but differs in several important aspects. First, the target rank $r$ and the sketch size $r+p$ are fixed \emph{a priori} and remain constant throughout the iterations. Second, the method stores only the last $k_{\text{trunc}}$ full basis tensors, while retaining only the sketches of the older basis tensors, since these are sufficient for reconstructing the final solution.

\paragraph{Memory efficient variant.}
However, with our choice of parameters, the sketches have size $\tilde{r} \coloneqq r+p$, which may even exceed the storage required by the corresponding basis tensors. Consequently, in memory-constrained settings or when a large number of iterations is required, storing all sketches can still become prohibitively expensive. While restarting can alleviate this issue, it may substantially impair convergence. We avoid this by exploiting the fact that the leading columns of $\bPsi_k \in \R^{n_k \times \tilde{r}}$ themselves form a valid, smaller Nystr\"om sketching operator. Specifically, if the desired solution rank is $r_{\mathrm{sol}} \le r$, we fix a (typically small) reduced oversampling parameter $p_{\mathrm{sol}}$ and set $\tilde{r}_{\mathrm{sol}} \coloneqq r_{\mathrm{sol}} + p_{\mathrm{sol}}$, retaining only the leading $\tilde{r}_{\mathrm{sol}}$ columns of each $\bPsi_k$. Of these, only the leading $r_{\mathrm{sol}}$ directions correspond to the target-rank contraction used to build $\bF_{\bPsi k}^{(j)}$ (\Cref{def:mttkrp}), which by construction carries no oversampling buffer of its own; the remaining $p_{\mathrm{sol}}$ columns are retained solely to preserve the numerical stability of the Nystr\"om reconstruction in \texttt{MLN\_Recovery}, mirroring the role that $p$ plays for the full-resolution sketch. We can therefore downsample every stored sketched core and factor to their corresponding leading indices,
\begin{equation}
\label{eq:memory_efficient_truncation}
    \bTG_{\bPsi}^{(j)} \leftarrow \bTG_{\bPsi}^{(j)}(1{:}\tilde{r}_{\mathrm{sol}}, \dots, 1{:}\tilde{r}_{\mathrm{sol}}),
    \qquad
    \bF_{\bPsi k}^{(j)} \leftarrow \bF_{\bPsi k}^{(j)}(:,\, 1{:}r_{\mathrm{sol}}), \quad k = 1, \dots, d,
\end{equation}
as done at line~\ref{lin:save_memory} of~\Cref{alg:mln_tucker_sgmres}. Since this is mathematically equivalent to having sketched with the reduced operator $\bPsi_k^{(\tilde{r}_{\mathrm{sol}})}$ from the start, the operation reduces the memory footprint of each stored basis tensor from $\mathcal{O}(d\tilde{r}^d)$ to $\mathcal{O}(d\tilde{r}_{\mathrm{sol}}^d)$ without altering the mathematical behavior or convergence of the algorithm, only the resolution at which the discarded basis tensors are retained. In our experiments we set $p_{\mathrm{sol}} = r_{\mathrm{sol}}$, doubling the retained resolution relative to the target solution rank. Moreover, as discussed above, the sketching required to solve the least-squares problem arises naturally as a byproduct of the MLN sketching phase, and therefore does not require the construction of additional sketching operators. Finally, unlike RHOSVD-Tucker GMRES, the sketch is applied to the exact operator application rather than to a truncated approximation, resulting in a more faithful representation of the projected residual, and does not introduce a mismatch with the truncated basis that may, in principle, affect the accuracy of the solution.

\begin{figure}[t]
    \centering
    \begin{minipage}[t]{0.49\textwidth}
        \vspace{0pt}
        \begin{algorithm}[H]
        \small
        \caption{\texttt{MLN\_Sketch}}
            \label{alg:mln_sketch}

            \SetKwInOut{Input}{Input}
            \SetKwInOut{Output}{Output}
            \DontPrintSemicolon

            \Input{Tucker tensors $\{\bTX^{(j)} = \bTG^{(j)} \times_{k=1}^d \bU_k^{(j)}\}_{j=1}^N$, \\
                   sketches $\{\bPsi_k\}_{k=1}^d$}
            \vspace{0.36em}
            \Output{Cores $\{\bTG_{\Psi}^{(j)}\}_{j=1}^N$,\\
                    factors $\{\{\bF_{\bPsi k}^{(j)}\}_{k=1}^d\}_{j=1}^N$,\\
                    (optional) $\mdiag(\bTG_{\bPsi})$}
            \vspace{0.36em}
            \For{$j \leftarrow 1$ \KwTo $N$}{
            \vspace{0.36em}
                \For{$k \leftarrow 1$ \KwTo $d$}{
                    $\bU_{\bPsi k}^{(j)} \leftarrow \bPsi_k^{\top} \bU_k^{(j)}$\;
                    $\overline{\bU}_{\bPsi k}^{(j)} \leftarrow \text{1st}\; r_k \text{ rows of } \bU_{\bPsi k}^{(j)}$\;
                }
                \vspace{0.36em}
                $\bTG_{\Psi}^{(j)} \leftarrow \bTG^{(j)} \times_{k=1}^d \bU_{\bPsi k}^{(j)}$\;
                \vspace{0.36em}
                \For{$k \leftarrow 1$ \KwTo $d$}{
                    $\bT \leftarrow \texttt{MTTKRP}\big(\bTG^{(j)},\{\overline{\bU}_{\bPsi \ell}^{(j)}\}_{\ell}, k\big)$\
                    $\bF_{\bPsi k}^{(j)} \leftarrow \bU_k^{(j)} \bT$\;
                }
                \vspace{0.36em}
            }
            \vspace{0.36em}
            \Return $\{\bTG_{\Psi}^{(j)}\}_{j=1}^N, \{\{\bF_{\bPsi k}^{(j)}\}_{k=1}^d\}_{j=1}^N$, (optional) $\mdiag(\bTG_{\bPsi})$\;
        \end{algorithm}
    \end{minipage}%
    \hfill%
    \begin{minipage}[t]{0.48\textwidth}
        \small
        \vspace{0pt}
        \begin{algorithm}[H]
            \caption{\texttt{MLN\_Sum}}
            \label{alg:mln_sum}

            \SetKwInOut{Input}{Input}
            \SetKwInOut{Output}{Output}
            \DontPrintSemicolon

            \Input{Cores $\{\bTG_{\Psi}^{(j)}\}_{j=1}^N$, factors $\{\{\bF_{\bPsi k}^{(j)}\}_{k=1}^d\}_{j=1}^N$,
                   coeffs. $\{c_j\}_{j=1}^N$}
            \Output{Core $\bTG_{\Psi}$,
                    factors $\{\bF_{\bPsi k}\}_{k=1}^d$}

            $\bTG_{\Psi} \leftarrow \sum_{j=1}^N c_j \bTG_{\Psi}^{(j)}$\;
            \For{$k \leftarrow 1$ \KwTo $d$}{
                $\bF_{\bPsi k} \leftarrow \sum_{j=1}^N c_j \bF_{\bPsi k}^{(j)}$\;
            }
            \Return $\bTG_{\Psi}, \{\bF_{\bPsi k}\}_{k=1}^d$\;
            
        \end{algorithm}
        \small
        \begin{algorithm}[H]
            \caption{\texttt{MLN\_Recovery}}
            \label{alg:mln_recovery}

            \SetKwInOut{Input}{Input}
            \SetKwInOut{Output}{Output}
            \DontPrintSemicolon

            \Input{Core $\bTG_{\Psi}$, factors $\{\bF_{\bPsi k}\}_{k=1}^d$, \\
                   sketches $\{\bPsi_k\}_{k=1}^d$}
            \Output{Core $\widetilde{\bTG}$,
                    factors $\{\widetilde{\bU}_k\}_{k=1}^d$}

            \For{$k \leftarrow 1$ \KwTo $d$}{
                $\bQ_k, \bR_k \leftarrow \texttt{qr}(\bPsi_k^{\top} \bF_{\bPsi k})$ \tcp*{Econ QR}
                $\widetilde{\bU}_k \leftarrow \bF_{\bPsi k} \bR_k^{\dagger}$\;
            }
            $\widetilde{\bTG} \leftarrow \bTG_{\Psi} \times_{k=1}^d \bQ_k^{\top}$\;
            \Return $\widetilde{\bTG}, \{\widetilde{\bU}_k\}_{k=1}^d$\;
            
        \end{algorithm}
    \end{minipage}
\end{figure}

\begin{algorithm}[p]
    \small
    \caption{\texttt{MLN-Tucker sGMRES}}
    \label{alg:mln_tucker_sgmres}
    \SetKwInOut{Input}{Input}
    \SetKwInOut{Output}{Output}
    \DontPrintSemicolon

    \Input{Linear operator $\opL$; Tucker tensor $\bTB$; target rank $r$; tolerance \texttt{tol}; sketched correction factor $\eta$; max iterations $N_{\max}$; truncation window $k_{\text{trunc}}$; \texttt{save\_memory} flag, if \texttt{True} do an additional truncation of the sketches based on
    ${r}_{\mathrm{sol}}$, $p_{\mathrm{sol}}$, and $\tilde{r}_{\mathrm{sol}} = {r}_{\mathrm{sol}}+p_{\mathrm{sol}}$}
    \Output{Approximate solution in Tucker format $\widetilde{\bTX}$, such that $\|\opL(\widetilde{\bTX})- \bTB\|_{\fro}\leq \texttt{tol}$}

    $\beta \gets \|\bTB\|_{\fro}$\;
    $\bTV_1 \gets \bTB / \beta$\;
    Generate Gaussian matrices $\{\bPsi_k\}_{k=1}^{d}$, with $\bPsi_k\in\mathbb{R}^{n_k\times r}$\;
    $\left[\bTG_{\Psi}^{(1)}, \{\bF_{\bPsi k}^{(1)}\}_{k=1}^d, \bb_s\right] \gets \texttt{MLN\_Sketch}\big(\{\bTV_1\}, \{\bPsi_k\}_{k=1}^d\big)$\;\label{lin: mln_sketchB}
    $\beta_s \gets \|\bb_s\|_2$\;
    $\bM \gets [\,\,]$\tcp*{initialization for   $\bS\bA\bV_k$}
    \For{$i \gets 1$ \KwTo $N_{\max}$}{
        $\bTW \gets \bTV_i$\;
        $\left[\bTG_{\Psi, \bTW}, \{\bF_{\bPsi k, \bTW}\}_{k=1}^d, \bs_{\bTW}\right] \gets \texttt{MLN\_Sketch}\big(\opL(\bTW), \{\bPsi_k\}_{k=1}^d\big)$\;
        $\bM \gets [\bM, \bs_{\bTW}]$\;\label{lin: mln_diag}

        $\bTW \gets \texttt{MLN\_Recovery}\big(\bTG_{\Psi, \bTW}, \{\bF_{\bPsi k, \bTW}\}_{k=1}^d, \{\bPsi_k\}_{k=1}^d\big)$\;

        \For{$j \gets \max(1, i-k_{\text{trunc}}+1)$ \KwTo $i$}{
            $h_{j,i} \gets \langle \bTV_j, \bTW \rangle$\;
            
            $\left[\bTG_{\Psi, \bTW}, \{\bF_{\bPsi k, \bTW}\}_{k=1}^d\right] \gets \texttt{MLN\_Sum}\Big( \{\bTG_{\Psi, \bTW}, \bTG_{\Psi}^{(j)}\}, \big\{\{\bF_{\bPsi k, \bTW}\}_{k=1}^d, \{\bF_{\bPsi k}^{(j)}\}_{k=1}^d\big\}, \{1, -h_{j,i}\} \Big)$\;

        }
        
        \If{$\|\bTW\|_{\fro} \approx 0$}{
            \KwBreak\;
        }
        
        $\bTV_{i+1} \gets \bTW / \|\bTW\|_{\fro}$\;
        $\left[\bTG_{\Psi}^{(i+1)}, \{\bF_{\bPsi k}^{(i+1)}\}_{k=1}^d\right] \gets \texttt{MLN\_Sketch}\big(\{\bTV_{i+1}\}, \{\bPsi_k\}_{k=1}^d\big)$\;

        \If{$i > k_{\text{trunc}}$}{
            Delete $\bTV_{i-k_{\text{trunc}}}$ from memory \tcp*{keep local orthogonality window}
                \If{\texttt{save\_memory}}{
                $\bTG_{\Psi}^{(i-k_\text{trunc})} = \bTG_{\Psi}^{(i-k_\text{trunc})}(1:\tilde{r}_{\mathrm{sol}}, \dots, 1:\tilde{r}_{\mathrm{sol}})$ and $\bF_{\bPsi k}^{(i-k_{\text{trunc}})} = \bF_{\bPsi k}^{(i-k_\text{trunc})}(:, 1:r_{\mathrm{sol}})$ for every $k = 1, \dots, d$\;
                }
        }
        
    $\by_i \gets \arg\min_{\by}\, \left\| \bM \by - \bb_s \right\|_2 (+\lambda\|\by\|_2^2$ if regularization is needed) 

        $\mathrm{res}_i \gets \|\bM \by_i - \bb_s\|_2/\beta_s$  \tcp*{Sketched relative residual}

        \If{$\mathrm{res}_i < \eta \cdot \textup{\texttt{tol}}$}{
            \KwBreak\; 
        }
    }  
    \If{\texttt{save\_memory}}{
    \For{$j = i - k_\text{trunc} + 1, \dots, i$}{
    $\bTG_{\Psi}^{(j)} = \bTG_{\Psi}^{(j)}(1:\tilde{r}_{\mathrm{sol}}, \dots, 1:\tilde{r}_{\mathrm{sol}})$ and $\bF_{\bPsi k}^{(j)} = \bF_{\bPsi k}^{(j)}(:, 1:r_{\mathrm{sol}})$ for $k = 1, \dots, d$\;\label{lin:save_memory}
    }}
    
    $\left[\bTG_{\Psi}^{\text{sol}}, \{\bF_{\bPsi k}^{\text{sol}}\}_{k=1}^d\right] \gets \texttt{MLN\_Sum}\Big( \{\bTG_{\Psi}^{(m)}\}_{m=1}^i, \big\{\{\bF_{\bPsi k}^{(m)}\}_{k=1}^d\big\}_{m=1}^i, \by_i \Big)$\;
    $\widetilde{\bTX} \gets \texttt{MLN\_Recovery}\big(\bTG_{\Psi}^{\text{sol}}, \{\bF_{\bPsi k}^{\text{sol}}\}_{k=1}^d, \{\bPsi_k\}_{k=1}^d\big)$\;

    \Return $\widetilde{\bTX}$\;
\end{algorithm}

\paragraph{Memory comparison.}
We compare the memory required by the solvers, focusing on the dominant contributions: the stored Krylov basis tensors in Tucker format and, where applicable, their sketched representations. Auxiliary quantities are negligible in comparison and are omitted from the following estimates. For the RHOSVD-based methods, the rank $r_i$ of the $i$-th basis tensor is not fixed but is instead estimated adaptively by the effective-rank heuristic within \texttt{RoundSum} (\Cref{sec: RHOSVD-Tucker GMRES}), and can vary, typically growing, across iterations; we denote by $r_{\max} = \max_i r_i$ its worst-case value over the run.
\begin{itemize}
    \item \textbf{RHOSVD-Tucker GMRES}: stores all $k$ basis tensors $\Rightarrow$ $\mathcal{O}\bigl(\sum_{i=1}^k (dnr_i + r_i^d)\bigr) = \mathcal{O}(k(dnr_{\max} + r_{\max}^d))$.
    \item \textbf{RHOSVD-Tucker sGMRES}: stores all $k$ basis tensors $\Rightarrow$ $\mathcal{O}\bigl(\sum_{i=1}^k (dnr_i + r_i^d)\bigr) = \mathcal{O}(k(dnr_{\max} + r_{\max}^d))$.
    \item \textbf{MLN-Tucker sGMRES}: stores only $k_{\mathrm{trunc}}$ basis tensors and their sketched representations $\Rightarrow$ $\mathcal{O}(k_{\mathrm{trunc}}(dnr + r^d) + k \cdot d\tilde{r}^d)$, where the second term accounts for the accumulated sketched cores and factors retained for the final solution assembly. Unlike the RHOSVD-based methods, here $r$ (and hence $\tilde r = r+p$) is fixed \emph{a priori} by construction rather than estimated, so this bound is exact rather than a worst case. Note that $\tilde{r}$ is always slightly larger than $r$, so this term can exceed the cost of storing the basis tensors themselves.
    \item \textbf{MLN-Tucker sGMRES} (memory efficient): stores only $k_{\mathrm{trunc}}$ basis tensors and \emph{downsampled} sketches (\Cref{eq:memory_efficient_truncation}) $\Rightarrow$ $\mathcal{O}(k_{\mathrm{trunc}}(dnr + r^d) + k \cdot d\, \tilde r_{\mathrm{sol}}^d)$.
\end{itemize}
Since $\tilde r_{\mathrm{sol}} \le \tilde r$, and typically $\tilde r_{\mathrm{sol}} \ll \tilde r$ whenever the solution is itself numerically low-rank, the memory efficient variant reduces the accumulated sketch cost substantially relative to the standard one; combined with $k_{\mathrm{trunc}} \ll k$, this yields the largest overall memory reduction among the four variants, and is the key advantage enabling MLN-Tucker sGMRES to reach higher Krylov-basis ranks or run for more iterations under a fixed memory budget. In contrast, because $r_{\max}$ for the RHOSVD-based methods is only known a posteriori and tends to grow with the iteration count, their memory footprint is comparatively harder to predict or bound in advance. 

\subsection{Preconditioning}\label{sec:Preconditionedlaplacian}

Assume that $\opM$ is an operator among tensor spaces of order and size compatible with $\opL$. Applying $\opM$ to a Krylov basis tensor $\bTV_k$ causes rank growth in the multilinear rank of the resulting tensor meaning we have to do some form of truncation $\mathscr{T}$ to keep computational tractability. Therefore the preconditioned state at step $k$ is actually given by $\bTZ_k = \mathscr{T}(\opM(\bTV_k))$ where the projection defined by $\mathscr{T}$ depends on the input tensor $\opM(\bTV_k)$. As such it is clear that the preconditioning operator $\opM_k(\cdot) := \mathscr{T}(\opM(\cdot))$ is nonlinear and changes in every iteration, thus we use the flexible GMRES (FGMRES) procedure~\cite{Saa93}. The only two variations from~\Cref{alg:tucker_gmres,alg:tucker_sgmres} are that we now save the preconditioned $\bTZ_k =  \opM_k(\bTV_k)$ before applying $\opL$ directly, and update $\widetilde{\bTX}$ directly within the preconditioned space as $\widetilde{\bTX} \gets \texttt{RoundSum}\left(\sum_{i=1}^k y_i \bTZ_i\right)$, thus bypassing the non-linear issues from the truncation step. Similarly in~\Cref{alg:mln_tucker_sgmres} we store and operate only with the core tensors $\bTG_{\Psi, \bTZ}^{(i)}$ and the factor matrices $\{\bF_{\bPsi k, \bTZ}^{(i)}\}_{k=1}^d$. A well-known preconditioning operator for PDEs results from an approximation of the inverse of the Laplacian operator, defined as
\begin{equation}
\label{eq:Lap}
    \bDelta_d=\bDelta_1\otimes \bI\otimes\cdots\otimes \bI +\bI\otimes\bDelta_1\otimes \cdots\otimes\bI + \cdots +
    \bI\otimes \cdots\otimes\bI\otimes\bDelta_1
\end{equation}
where $\bI$ is the identity matrix, and $\bDelta_1$ is the $(n+1\times n+1)$ discretized 1-dimensional Laplacian over a grid of step size $h = 1/n$, i.e., $\bDelta_1 = (n+1)^{2}{\rm tridiag}(-1, 2, -1)$.
In~\cite{HacK06}, the authors prove that the approximation of the inverse of the $d$-dimensional Laplacian can be expressed as
\begin{equation}
    \bDelta^{-1}_{d} = \sum_{h=-q}^{q} c_h \exp(t_h\bDelta_1) \otimes\cdots\otimes\exp(t_h\bDelta_1)
\end{equation}
where $c_h = \eta t_h$, $t_h = \exp(h\eta)$, $\eta = \pi/\sqrt{q}$. The action of this operator can be efficiently computed via spectral decomposition. 
Let $\bV$ and $\bLambda$ denote the eigenvectors and eigenvalues of $\bDelta_1$. Then
\begin{equation}
    \exp(t_h\bDelta_1) \otimes \cdots \otimes \exp(t_h\bDelta_1)
    = \opV \exp\bigl(t_h(\bLambda \oplus \cdots \oplus \bLambda)\bigr)\opV^{-1},
\end{equation}
where $\opV = \bV \otimes \cdots \otimes \bV$, see~\cite{Hig08,Dol13}. For a Tucker tensor $\bTX = \bTC \times_{k=1}^d \bU_k$, the transformations $\opV^{-1}\bTX$ and $\opV\bTX$ reduce to applying the discrete sine transform of type I (DST-I) and its inverse to the factor matrices $\bU_k$ yielding an efficient implementation of the preconditioner.


\section{Numerical Experiments}
\label{sec: Results}

The reported numerical experiments in this section were performed with computational resources provided by Advanced Research Computing (ARC) at Virginia Tech (see acknowledgments) using a system with 128\,GB of RAM running Linux (5.14.0-687.26.1.el9\_8.x86\_64) with MATLAB 25.2.0.2998904 (R2025b). The source codes, data, and results of the numerical experiments reported in
this section are available at~\cite{supBucIPetal26}.

We consider three test problems and compare the proposed RHOSVD-Tucker sGMRES and MLN-Tucker sGMRES methods with standard low-rank solvers, including RHOSVD-Tucker GMRES and, in the symmetric case, the Subspace Gradient Descent method in Tucker format (TK-SS-SD)\cite{IanPS26}.
Firstly we test on the preconditioned heat equation using the Laplacian preconditioner described in Section~\ref{sec:Preconditionedlaplacian} which typically favors methods that exploit symmetry, such as conjugate gradient variants, due to their short recurrences and low memory requirements, while GMRES is disadvantaged by its long recurrences. 
Secondly we consider a 3D convection-diffusion equation wherein we break the symmetry of the Laplacian operator and introduce a convection term, increasing the sensitivity of Krylov methods to orthogonalization. We do not, however, precondition this example, so that we may test the effect of a larger number of iterations on the solvers where preconditioners may not be available.
Finally, we consider an inverse problem arising from image deblurring. By constructing a low-rank approximation to the right hand side we are able to compute a solution to the inverse problem which usually is intractable. While this truncation introduces an additional source of error, it also acts as a regularizer which we demonstrate.

\subsection{Laplacian}
\label{sec:laplacian}

The first problem we consider is a classical Poisson equation of the form $-\Delta u = f$ where $u$ is a function defined on a 3D domain, i.e., $u \coloneqq u(x,y,z)$. Discretizing it with a central finite-difference scheme, we obtain a tensor equation of the form~\eqref{eq:AXB:teneq} where the operator $\mathcal{L}$ corresponds to $\bDelta_3$ defined in~\Cref{eq:Lap} with $d=3$. The right-hand is obtained by discretizing $f$ and is represented in separable form as $f = \bv_1 \otimes \bv_2 \otimes \bv_3$, where $\bv_1$ is the normalized vector of all ones, and $\bv_2$, $\bv_3$ are the first canonical basis vectors. The chosen grid size is $n=1000$. We focus on preconditioned iterations and use exponential sum preconditioners tailored to the structure of~\eqref{eq:Lap}, described in~\Cref{sec:Preconditionedlaplacian}. In~\Cref{tab:prec_laplacian_summary}, we report the solvers' convergence at different accuracy values $10^{-3}, 10^{-5}, 10^{-7}$. For the three tolerances, we set the TK-SS-SD maximum ranks to $10$, $20$, and $30$ with truncation tolerance of $10^{-12}$. For the sketched solvers, we set the Nystr\"om sketching ranks to $10$, $25$, and $40$ with standard oversampling parameters of $p = 10, 25,$ and $40$, respectively, yielding $\tilde{r} = 20, 50, 80$. For the memory-efficient MLN variant, we restrict $\tilde{r}_{\mathrm{sol}} = 20,40,60$ to be made up of $r_{\mathrm{sol}} = p_{\mathrm{sol}} = 10, 20,  30$. Additionally, we choose the GMRES truncation index of $k_\text{trunc}= 1$, the sketching $p$-rank equal to $5$, the sketching tolerance multiplier $\eta_{\text{sketch}}=0.4$.

In the low and medium accuracy setting, the convergence time of the GMRES variants is comparable to the TK-SS-SD. In the high accuracy case, the GMRES algorithms outperform the SS-SD solver, which was originally designed for a low-to-medium accuracy setting. The small accuracy drop in the memory efficient version demonstrates that sometimes you have to be careful choosing what rank to truncate the stored sketches to for reconstruction. Here we have chosen too small a reconstruction rank to demonstrate that such an accuracy drop can happen. Nonetheless, all methods in the finer tolerance scheme perform very well.

\begin{table}[t]
    \centering
     \caption{Laplacian example: CPU time in seconds, algorithm iteration count, and true relative residual at convergence for differing tolerance levels.}
  \tikzexternalenable%
  \tikzsetnextfilename{prec_laplacian_summary}%
  \begin{tikzpicture}[font = \plotfontsize\normalfont]
  
  \pgfplotstableread[col sep=tab]{graphics/data/laplacian_summary.dat}\tablelong
  
  \pgfplotstablenew[
    columns={Method, SortOrder, TimeIter3, Res3, TimeIter5, Res5, TimeIter7, Res7},
    create on use/SortOrder/.style={
      create col/assign/.code={
        \ifcase\pgfplotstablerow
          \def\myorder{1} 
        \or
          \def\myorder{3} 
        \or
          \def\myorder{4} 
        \or
          \def\myorder{5} 
        \or
          \def\myorder{2} 
        \fi
        \pgfkeyslet{/pgfplots/table/create col/next content}{\myorder}
      }
    },
    create on use/Method/.style={
      create col/assign/.code={
        \ifnum\pgfplotstablerow=0
          \def\mymethod{TK-SS-SD}
        \else
          \pgfplotstablegetelem{\pgfplotstablerow}{Method}\of{\tablelong}
          \let\mymethod\pgfplotsretval
        \fi
        \pgfkeyslet{/pgfplots/table/create col/next content}{\mymethod}
      }
    },
    create on use/TimeIter3/.style={
      create col/assign/.code={
        \pgfmathtruncatemacro{\idx}{\pgfplotstablerow}
        \pgfplotstablegetelem{\idx}{Time(s)}\of{\tablelong}
        \pgfmathprintnumberto[fixed, precision=4]{\pgfplotsretval}{\timeval}
        \pgfplotstablegetelem{\idx}{Iterations}\of{\tablelong}
        \edef\combined{\timeval\ (\pgfplotsretval)}
        \pgfkeyslet{/pgfplots/table/create col/next content}{\combined}
      }
    },
    create on use/Res3/.style={
      create col/assign/.code={
        \pgfmathtruncatemacro{\idx}{\pgfplotstablerow}
        \pgfplotstablegetelem{\idx}{Final_TrueRes}\of{\tablelong}
        \pgfkeyslet{/pgfplots/table/create col/next content}{\pgfplotsretval}
      }
    }, 
    create on use/TimeIter5/.style={
      create col/assign/.code={
        \pgfmathtruncatemacro{\idx}{\pgfplotstablerow + 5} 
        \pgfplotstablegetelem{\idx}{Time(s)}\of{\tablelong}
        \pgfmathprintnumberto[fixed, precision=4]{\pgfplotsretval}{\timeval}
        \pgfplotstablegetelem{\idx}{Iterations}\of{\tablelong}
        \edef\combined{\timeval\ (\pgfplotsretval)}
        \pgfkeyslet{/pgfplots/table/create col/next content}{\combined}
      }
    },
    create on use/Res5/.style={
      create col/assign/.code={
        \pgfmathtruncatemacro{\idx}{\pgfplotstablerow + 5}
        \pgfplotstablegetelem{\idx}{Final_TrueRes}\of{\tablelong}
        \pgfkeyslet{/pgfplots/table/create col/next content}{\pgfplotsretval}
      }
    },
    create on use/TimeIter7/.style={
      create col/assign/.code={
        \pgfmathtruncatemacro{\idx}{\pgfplotstablerow + 10} 
        \pgfplotstablegetelem{\idx}{Time(s)}\of{\tablelong}
        \pgfmathprintnumberto[fixed, precision=4]{\pgfplotsretval}{\timeval}
        \pgfplotstablegetelem{\idx}{Iterations}\of{\tablelong}
        \edef\combined{\timeval\ (\pgfplotsretval)}
        \pgfkeyslet{/pgfplots/table/create col/next content}{\combined}
      }
    },
    create on use/Res7/.style={
      create col/assign/.code={
        \pgfmathtruncatemacro{\idx}{\pgfplotstablerow + 10}
        \pgfplotstablegetelem{\idx}{Final_TrueRes}\of{\tablelong}
        \pgfkeyslet{/pgfplots/table/create col/next content}{\pgfplotsretval}
      }
    }
  ]{5}\tablewide

  \node (timingtable) at (0,0) {
    \resizebox{\textwidth}{!}{%
    \begingroup
    \setlength{\tabcolsep}{3pt} 
    \footnotesize 
    \pgfplotstabletypeset[
      sort=true,
      sort key=SortOrder,
      every head row/.style={
        before row={
          \toprule
          & \multicolumn{2}{c}{\textbf{Tol:} $\mathbf{10^{-3}}$} 
          & \multicolumn{2}{c}{\textbf{Tol:} $\mathbf{10^{-5}}$} 
          & \multicolumn{2}{c}{\textbf{Tol:} $\mathbf{10^{-7}}$} \\
          \cmidrule(lr){2-3} \cmidrule(lr){4-5} \cmidrule(lr){6-7}
        },
        after row=\midrule
      },
      every last row/.style={after row=\bottomrule},
      columns={Method, TimeIter3, Res3, TimeIter5, Res5, TimeIter7, Res7},
      columns/Method/.style={
        string type, 
        column type=l, 
        column name=\textbf{Method}
      },
      columns/TimeIter3/.style={string type, column type=c, column name={\textbf{Time (Iter)}}},
      columns/Res3/.style={sci, sci zerofill, precision=2, column type=c, column name={\textbf{True Res.}}},
      columns/TimeIter5/.style={string type, column type=c, column name={\textbf{Time (Iter)}}},
      columns/Res5/.style={sci, sci zerofill, precision=2, column type=c, column name={\textbf{True Res.}}},
      columns/TimeIter7/.style={string type, column type=c, column name={\textbf{Time (Iter)}}},
      columns/Res7/.style={sci, sci zerofill, precision=2, column type=c, column name={\textbf{True Res.}}},
    ]{\tablewide}
    \endgroup
    }%
  };
\end{tikzpicture}%
  \tikzexternaldisable%

    \label{tab:prec_laplacian_summary}
\end{table}

\subsection{Convection-diffusion}
\label{sec:convection}

\begin{figure}[t]
    \centering
    \begin{subfigure}[b]{0.48\textwidth}
        \centering
  \tikzexternalenable%
  \tikzsetnextfilename{no_prec_convection_true_residual}%
  \begin{tikzpicture}[font = \normalfont]

  \pgfplotstableread[col sep=tab]{graphics/data/convection_MLNTuckersGMRES_trueres_k2.dat}\tableMLNsGMRES
  \pgfplotstableread[col sep=tab]{graphics/data/convection_RHOSVDTuckersGMRES_trueres_k2.dat}\tablesGMRES
  \pgfplotstableread[col sep=tab]{graphics/data/convection_RHOSVDTuckerGMRES_trueres_k2.dat}\tableGMRES
  \pgfplotstableread[col sep=tab]{graphics/data/convection_MLNTuckersGMRESMemoryEfficient_trueres_k2.dat}\tableMLNsGMRESMemoryEfficient
  
  \begin{axis}[%
    width                  = \textwidth,
    height                 = 5.5cm,
    ymode                  = log,
    xmin                   = 0,
    xmax                   = 110,
    ymin                   = 1e-5,
    ymax                   = 1e+0,
    ytick                  = {1e+0, 1e-2, 1e-4},
    xlabel                 = {Iteration},
    ylabel                 = {True Residual},
    grid                   = major,
    legend pos             = north east,
    legend style           = {font=\scriptsize, inner sep=2pt, legend cell align=left}
  ]   
  
    \addplot[GMRES, mark repeat = 8] 
      table[x = Iteration, y = True_Residual] {\tableGMRES};
    \addlegendentry{RHOSVD GMRES} 

    \addplot[sGMRES, mark repeat = 5] 
      table[x = Iteration, y = True_Residual] {\tablesGMRES};
    \addlegendentry{RHOSVD sGMRES} 

    \addplot[MLNsGMRES, mark repeat = 7] 
      table[x = Iteration, y = True_Residual] {\tableMLNsGMRES};
    \addlegendentry{MLN sGMRES} 

    \addplot[MLNsGMRESMemoryEfficient, mark repeat = 6] 
      table[x = Iteration, y = True_Residual] {\tableMLNsGMRESMemoryEfficient};
    \addlegendentry{\begin{tabular}{@{}l@{}}MLN sGMRES \\ (Memory Efficient)\end{tabular}} 
  \end{axis}
\end{tikzpicture}%
  \tikzexternaldisable%

        \caption{True residual comparison}
        \label{fig:no_prec_convection_true_residual}
    \end{subfigure}
    \hfill
    \begin{subfigure}[b]{0.48\textwidth}
        \centering
  \tikzexternalenable%
  \tikzsetnextfilename{no_prec_convection_ktrunc}%
  \begin{tikzpicture}[font = \normalfont]

  \pgfplotstableread[col sep=tab]{graphics/data/convection_RHOSVDTuckersGMRES_trueres_ktrunc1.dat}\tableKone
  \pgfplotstableread[col sep=tab]{graphics/data/convection_RHOSVDTuckersGMRES_trueres_ktrunc2.dat}\tableKtwo
  \pgfplotstableread[col sep=tab]{graphics/data/convection_RHOSVDTuckersGMRES_trueres_ktrunc3.dat}\tableKthree

  \begin{axis}[%
    width                  = \textwidth,
    height                 = 5.5cm,
    ymode                  = log,
    xmin                   = 0,
    xmax                   = 110,
    ymin                   = 1e-6,
    ymax                   = 1e+0,
    ytick                  = {1e+0, 1e-2, 1e-4, 1e-6},
    xlabel                 = {Iteration},
    ylabel                 = {True Residual},
    grid                   = major,
    legend pos             = north east
  ]
  
    \addplot[sGMRES, mark repeat = 5, color=matlabred, mark=square*] 
      table[x = Iteration, y = True_Residual] {\tableKone};
    \addlegendentry{$\texttt{k\_trunc} = 1$} 

    \addplot[sGMRES, mark repeat = 6, color=matlablightblue, mark=*] 
      table[x = Iteration, y = True_Residual] {\tableKtwo};
    \addlegendentry{$\texttt{k\_trunc} = 2$} 

    \addplot[sGMRES, mark repeat = 7, color=matlabpurple, mark=triangle*] 
      table[x = Iteration, y = True_Residual] {\tableKthree};
    \addlegendentry{$\texttt{k\_trunc} = 3$} 
    
  \end{axis}
\end{tikzpicture}%
  \tikzexternaldisable%

        \caption{Effect of $k_{\mathrm{trunc}}$ on RHOSVD sGMRES}
        \label{fig:no_prec_convection_ktrunc}
    \end{subfigure}
    
    \caption{Convergence profiles for the unpreconditioned 3D convection-diffusion problem. (a) True residual histories comparing the proposed MLN and RHOSVD sketched solvers against the full-orthogonalization baseline. (b) Influence of the partial orthogonalization window $k_{\mathrm{trunc}}$ on RHOSVD-Tucker sGMRES, demonstrating residual stagnation for $k_{\mathrm{trunc}}=1$ and robust convergence for $k_{\mathrm{trunc}} \ge 2$.}
    \label{fig:combined_convection_metrics}
\end{figure}
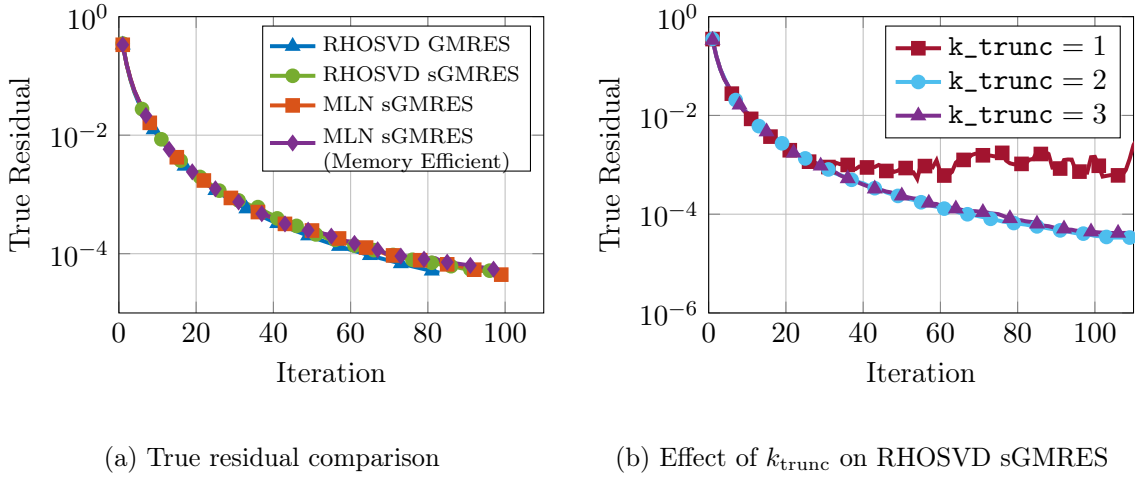

We consider a 3D convection-diffusion equation of the form $-\nabla \cdot (\kappa\nabla u) + \bv \cdot \nabla u = f$ defined on the unit cube domain $[0, 1]^3$.  The non-symmetric convection term yields a non-self-adjoint operator, altering the spectral properties of the system and increasing the sensitivity of Krylov methods to the orthogonalization strategy, making it a good benchmark in tandem with~\Cref{sec:Preconditionedlaplacian}. Using the Kronecker sum structure, the discrete operator $\opL$ in~\eqref{eq:AXB:teneq} takes the separable form
\begin{equation}
    \opL = \bL_1 \otimes \bI \otimes \bI + \bI \otimes \bL_1 \otimes \bI + \bI \otimes \bI \otimes \bL_1,
\end{equation}
where $\bL_1 = \kappa \bDelta_1 + \omega \bNabla_1 \in \R^{n \times n}$ is the 1D discrete convection-diffusion operator, with $\omega$ the uniform scalar constant velocity in any single direction from $\bv$. Here, $\bDelta_1 = \frac{1}{h^2} \mathrm{tridiag}(-1, 2, -1)$ represents the 1D discrete negative Laplacian, and $\bNabla_1$ is the backward difference matrix resulting from an upwind convection scheme, represented by subdiagonals $[-1, 0]$ with values $[-1/h, 1/h]$. We set the coefficients to $\kappa = 10^{-2}$ and $\omega = 5 \times 10^{-2}$. The right-hand side is a rank-1 tensor constructed as $\bTB = \ba_1 \otimes \bb_1 \otimes \bc_1$, where $\ba_1 = \be / \|\be\|_2$ (with $\be$ being the vector of all ones), and $\bb_1 = \bc_1 = \be_1 / \|\be_1\|_2$ are normalized standard basis vectors. The chosen grid size is $n=1000$. For this benchmark, we target a relative residual tolerance of $5 \times 10^{-6}$ over a maximum of $120$ iterations. We set the Nystr\"om target rank to $60$ with a standard oversampling $p = 180$. For the memory-efficient MLN variant, we use $r_{\mathrm{sol}} = p_{\mathrm{sol}}= 25$ yielding $\tilde{r}_{\mathrm{sol}} = 50$, a significant reduction in stored sketch dimension for the cores. Additionally, we choose a Khatri-Rao sketching $p$-rank equal to $5$ and a sketching tolerance multiplier $\eta_{\text{sketch}}=0.3$. For the main convergence comparison~\Cref{fig:no_prec_convection_true_residual}, the GMRES truncation index is fixed at $k_\text{trunc}= 2$; however, to analyze the sensitivity of the short recurrences to the non-symmetric convection term~\Cref{fig:no_prec_convection_ktrunc}, we also evaluate the RHOSVD-Tucker sGMRES method across an orthogonalization window of $k_\text{trunc} \in \{1, 2, 3\}$.

Here we compare our proposed methods against standard RHOSVD-Tucker GMRES on the unpreconditioned system, demonstrating not only true residual convergence of the proposed methods with times and iteration count, but also the effect of increasing $k_\text{trunc}$ in~\Cref{alg:mln_tucker_sgmres,alg:tucker_sgmres} in the short recurrence. Our results in~\Cref{fig:no_prec_convection_ktrunc} empirically demonstrate that using $k_\text{trunc} = 2$ is a good choice, and that under this choice the methods' true residuals converge in~\Cref{fig:no_prec_convection_true_residual}. However, the higher iteration count here compared to the preconditioned Laplacian in~\Cref{sec:Preconditionedlaplacian} causes some issues. Namely, the faster accumulation of rank in RHOSVD-Tucker GMRES causes a slower run time compared to its sketched counterparts as seen in~\Cref{tab:no_prec_convection_summary}. Although RHOSVD-Tucker sGMRES outperforms in solve time, it does not show true residual convergence, namely premature convergence, so $\eta = 0.3$ may not be strong enough here. The extra 3 iterations here would yield similar timings therefore to the MLN methods that do reach desired tolerance. Importantly, in the large iteration scheme we have successfully saved the memory in the memory efficient version by only storing the reconstruction sketches according to $\tilde{r}_{\mathrm{sol}} = 50$ rather than $\tilde{r} = 240$, a significant improvement. The memory efficient method nicely follows the true residual of the full MLN method closely, converging in the same number of iterations, and yielding a similar solve time. 

\begin{table}[t]
    \centering
    \caption{Performance summary for the unpreconditioned 3D convection-diffusion problem, comparing execution time, final true residual with target tolerance $5 \times 10^{-5}$, and iteration count across the proposed and baseline Tucker-GMRES variants.}
  \tikzexternalenable%
  \tikzsetnextfilename{no_prec_convection_summary}%
  \begin{tikzpicture}[font = \plotfontsize\normalfont]
  
  \pgfplotstableread[col sep=tab]{graphics/data/convection_summary.dat}\tableconvectionsummary

  \node (timingtable) at (0,0) {
    \pgfplotstabletypeset[
      every head row/.style={
        before row={\toprule},
        after row=\midrule
      },
      every last row/.style={after row=\bottomrule},
      columns={Method, k_trunc, {Time(s)}, Final_TrueRes, Iterations},
      columns/Method/.style={
        string type, 
        column type=l, 
        column name=\textbf{Method}
      },
      columns/k_trunc/.style={
        column type=c, 
        column name={\textbf{$\texttt{k\_trunc}$}},
        assign cell content/.code={
          \pgfplotstablegetelem{\pgfplotstablerow}{Method}\of{\tableconvectionsummary}%
          \edef\currentmethod{\pgfplotsretval}%
          \def\targetmethodA{RHOSVD-Tucker GMRES}%
          \def\targetmethodB{TuckerAxby}%
          \ifx\currentmethod\targetmethodA
            \pgfkeyssetvalue{/pgfplots/table/@cell content}{-}%
          \else\ifx\currentmethod\targetmethodB
            \pgfkeyssetvalue{/pgfplots/table/@cell content}{-}%
          \else
            \pgfkeyssetvalue{/pgfplots/table/@cell content}{##1}%
          \fi\fi
        }
      },
      columns/{Time(s)}/.style={
        fixed, 
        fixed zerofill, 
        precision=2,
        column type=c, 
        column name={\textbf{Solve Time (s)}}
      },
      columns/Final_TrueRes/.style={
        sci, 
        sci zerofill, 
        precision=2, 
        column type=c, 
        column name={\textbf{True Res.}}
      },
      columns/Iterations/.style={
        column type=c, 
        column name={\textbf{Iter.}}
      },
      sort=false,
      sort key={Method},
      sort cmp={string <}
    ]{\tableconvectionsummary}
  };
y\end{tikzpicture}%
  \tikzexternaldisable%

    \label{tab:no_prec_convection_summary}
\end{table}

\subsection{Image Deblurring}
\label{sec:image_deblurring}

\begin{figure}[t]
\centering
\begin{subfigure}[t]{0.3\textwidth}
\centering
\includegraphics[width=\textwidth]{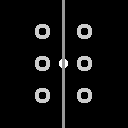}
\caption{}
\label{fig:deblur_true}
\end{subfigure}
\hfill
\begin{subfigure}[t]{0.3\textwidth}
\centering
\includegraphics[width=\textwidth]{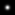}
\caption{}
\label{fig:deblur_psf}
\end{subfigure}
\hfill
\begin{subfigure}[t]{0.3\textwidth}
\centering
\includegraphics[width=\textwidth]{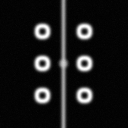}
\caption{}
\label{fig:deblur_blurred}
\end{subfigure}
\caption{Image deblurring problem.
(a)~Ground truth at time instance $n_d = 31$ ($128 \times 128$ pixels).(b)~PSF ($15 \times 15$ pixels)
(c)~Blurred and noisy image.}
\label{fig:deblur_setup}
\end{figure}
Our third example addresses a large-scale inverse problem arising from three-dimensional image deblurring.  Unlike the PDE-driven problems of the previous sections, inverse problems are characterized by the presence of noise in the observed data and the severe ill-conditioning of the forward operator $\opL$, whose singular values decay gradually to zero~\cite{Han10, HanNO06}.  As a consequence, naively solving~\eqref{eq:AXB:teneq} strongly amplifies the noise, and classical iterative methods exhibit \emph{semi-convergence}: the reconstruction error first decreases as dominant spectral components are captured, then increases as high-frequency noise contaminates the approximation~\cite{HanNO06, Han10}.  Regularization is therefore essential to obtain a stable reconstruction of the unknown $\bTX_{\mathrm{true}}$.
To assess reconstruction quality when the true solution is available, we use the relative reconstruction error (RRE), defined as $\lVert \bTX_k - \bTX_{\mathrm{true}} \rVert_{\fro} / \lVert
 \bTX_{\mathrm{true}} \rVert_{\fro}$ where $\bTX_k$ denotes the reconstructed solution at iteration $k$.            
 
When $n_d = 1$, i.e., a single two-dimensional image, the blurring operator can be represented as a matrix $\bA \in \R^{N \times N}$ with $N = n_1 n_2$, and hybrid projection methods that embed regularization within the Krylov iteration. Methods such as hybrid GMRES~\cite{GazNR15} have proven effective~\cite{Han10, CalGR99}.  However, when $n_d \gg 1$ (e.g., a video sequence), the vectorized system has dimension $N = n_1 n_2 n_d$, and storing the full Arnoldi basis and the dense right-hand side becomes prohibitively expensive. Working entirely in the Tucker format offers a path to tractability: it reduces storage from $\mathcal{O}(N)$ to $\mathcal{O}(\sum_k n_k r_k + \prod_k r_k)$, and, as we demonstrate below, low-rank compression of the data acts as an implicit form of regularization.
 
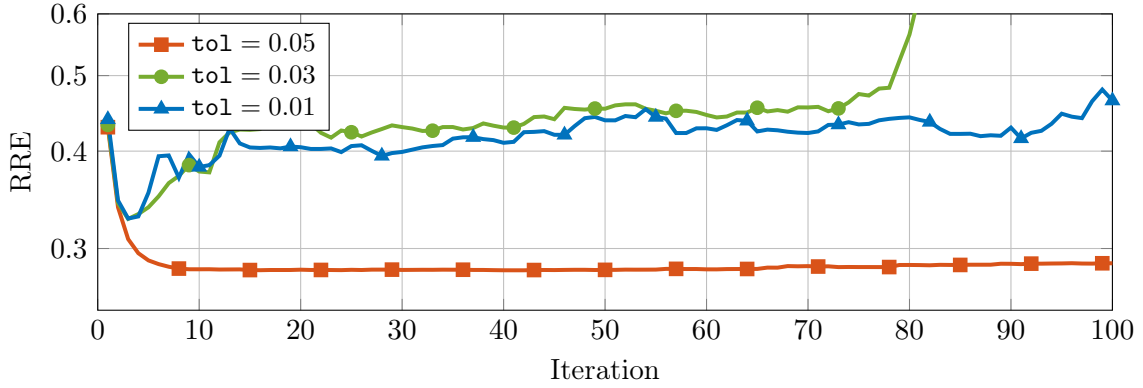
\begin{figure}[t]
    \centering
  \tikzexternalenable%
  \tikzsetnextfilename{deblurring_residuals_comp}%
  \begin{tikzpicture}

    \pgfplotstableread[col sep=tab]{graphics/data/deblurring_comp_MLNTuckersGMRES_history_k2_tol5em02.dat}\tableMLNtolA
    \pgfplotstableread[col sep=tab]{graphics/data/deblurring_comp_MLNTuckersGMRES_history_k2_tol3em02.dat}\tableMLNtolB
    \pgfplotstableread[col sep=tab]{graphics/data/deblurring_comp_MLNTuckersGMRES_history_k2_tol1em02.dat}\tableMLNtolC

    \begin{axis}[
        width             = \textwidth,
        height            = 5.5cm,
        ymode             = log,
        ytick             ={0.3, 0.4, 0.5, 0.6},
        xmin              = 0,
        xmax              = 100,
        ymin              = 0.25,
        ymax              = 0.6,
        minor y tick num  = 5,
        log ticks with fixed point,
        xlabel            = {Iteration},
        ylabel            = {RRE},
        grid              = major,
        legend pos        = north west,
        legend style      ={font=\small},
    ]
    
    \addplot[MLNsGMRES, mark repeat=7, color=matlaborange, mark=square*] 
        table[x=Iteration, y=RRE] {\tableMLNtolA};
    \addlegendentry{$\texttt{tol}=0.05$}
        
    \addplot[MLNsGMRES, mark repeat=8, color=matlabgreen, mark=*] 
        table[x=Iteration, y=RRE] {\tableMLNtolB};
    \addlegendentry{$\texttt{tol}=0.03$}
        
    \addplot[MLNsGMRES, mark repeat=9, color=matlabblue, mark=triangle*] 
        table[x=Iteration, y=RRE] {\tableMLNtolC};
    \addlegendentry{$\texttt{tol}=0.01$}
        
    \end{axis}
\end{tikzpicture}%
  \tikzexternaldisable%

    \caption{RRE history showing the effect of varying the RHS compression tolerance on the MLN-Tucker sGMRES method.}
    \label{fig:deblurring_res_comp}
\end{figure}

We consider a ground truth tensor $\bTX_{\mathrm{true}}\in\R^{128\times
128\times 64}$ given by a synthetic \emph{hollow bars} phantom inspired by
the EPFL DeconvolutionLab benchmark~\cite{KirASetal13, EpfND}.  It contains six hollow cylindrical bars arranged on a $3\times 2$ grid with outer radius~8 and inner radius~4 pixels.  Each bar has an axial
intensity taper so that it smoothly vanishes near the top and bottom of the temporal slices.  At the volume center the phantom contains a small spherical bead of radius~5 and a thin planar membrane at the mid-$y$ plane.  All intensities are normalized to $[0,1]$.  The true image at time instance $31$ is shown in~\Cref{fig:deblur_true}. The widefield point-spread function (PSF), displayed in~\Cref{fig:deblur_psf}, is modeled as a separable
product of three one-dimensional Gaussians,
\begin{align*}
h(i,j,k) &= h_1(i) \, h_2(j) \, h_3(k) \\
\text{where:} \quad h_1(t) = h_2(t) &= \frac{1}{Z_{xy}} \exp\left(-\frac{t^2}{2\sigma_{xy}^2}\right) \\
h_3(t) &= \frac{1}{Z_z} \exp\left(-\frac{t^2}{2\sigma_z^2}\right)
\end{align*}
with $\sigma_{xy}=1.5$ pixels (lateral), $\sigma_z=3.5$ pixels (axial),
and $Z_{xy}, Z_z$ are normalization constants ensuring each kernel sums to
one.
Each 1-D kernel is embedded into an $n_k\times n_k$ Toeplitz-like
blurring matrix $\bA_k$, $k =1, 2, 3$, with reflective (Neumann) boundary
conditions.  The blurred, simulated tensor is then obtained via successive mode-$k$ products,
 $ \bTY_{\rm true}
  \;=\;
 \bTX_{\mathrm{true}} \times_1 \bA_1 \times_2 \bA_2 \times_3 \bA_3$.
In vectorized form this reads
\begin{equation*}\label{eq:forward_vec}
  \mathrm{vec}(\bTY_{\rm true})
  \;=\;
\underbrace{(\bA_3\otimes\bA_2\otimes\bA_1)}_{\displaystyle \bA\;\in\;\R^{N\times N}}
  \;\mathrm{vec}(\bTX_{\mathrm{true}}).
\end{equation*}
where $N = n_1 n_2 n_3$.
The Kronecker structure of $\bA$ is exploited so that the matrix is never formed explicitly; its action on a vector (or, equivalently, on a Tucker tensor) is evaluated through sequential mode-$k$ multiplications with $\bA_1$, $\bA_2$, and $\bA_3$. Gaussian noise is added at a controlled level such that
\begin{equation*}\label{eq:noise}
  \bTY
  \;=\;
  \bTY_{\rm true} + \delta\,\frac{\|\bTY_{\rm true}\|_{\fro}}{\|\boldsymbol{\eta}\|_{\fro}}\,\boldsymbol{\eta}, \quad \boldsymbol{\eta}\sim\mathcal{N}(0,\bI),
\end{equation*}
where $\delta = 0.05$ is the noise level and $\bI$ denoting the identity covariance matrix. The blurred and noisy image at time frame 31 is shown in ~\Cref{fig:deblur_blurred}. 
 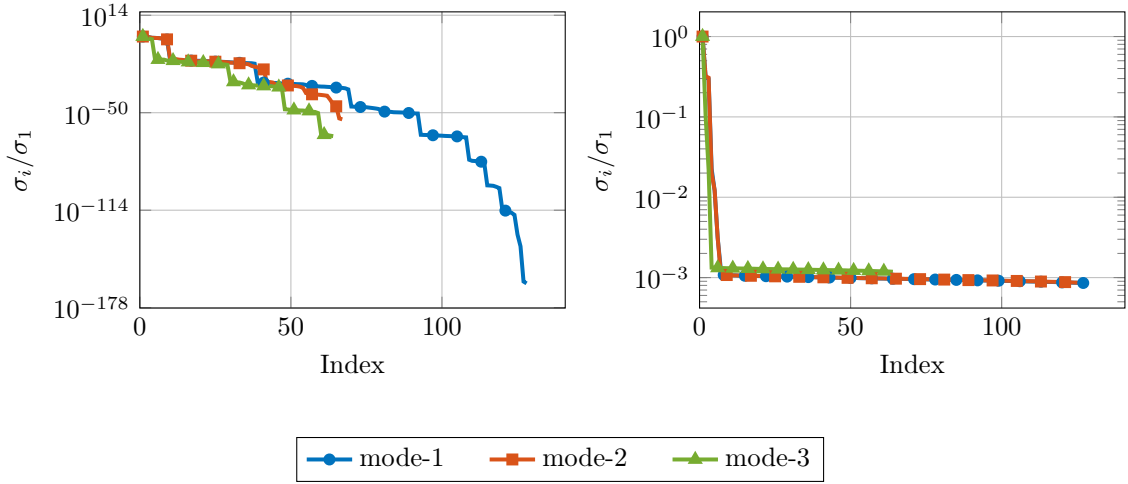
\begin{figure}[t]
      \centering
      \begin{subfigure}[t]{0.48\textwidth}
          \centering                %
  \tikzexternalenable%
  \tikzsetnextfilename{svd_unfoldings_xtrue}%
  \begin{tikzpicture}[font = \plotfontsize\normalfont]

  \pgfplotstableread[col sep=tab]{graphics/data/svd_unfoldings_Xtrue.dat}\tableXtrue

  \begin{axis}[%
    width            = \textwidth,
    height           = 5.5cm,
    ymode            = log,
    xmin             = 0,
    xlabel           = {Index},
    ylabel           = {$\sigma_i / \sigma_1$},
    grid             = major,
  ]

    \addplot[matlabblue, solid, line width=1.5pt,
             mark=*, mark options={solid}, mark repeat=8, mark size=1.5pt]
      table[x=Index, y=Mode1] {\tableXtrue};

    \addplot[matlaborange, solid, line width=1.5pt,
             mark=square*, mark options={solid}, mark repeat=8, mark size=1.5pt]
      table[x=Index, y=Mode2] {\tableXtrue};

    \addplot[matlabgreen, solid, line width=1.5pt,
             mark=triangle*, mark options={solid}, mark repeat=5, mark size=1.8pt]
      table[x=Index, y=Mode3] {\tableXtrue};

  \end{axis}
\end{tikzpicture}%
  \tikzexternaldisable%

          \label{fig:svd_unfolding_xtrue}                           \end{subfigure}          
      \hfill                                
      \begin{subfigure}[t]{0.48\textwidth}              
          \centering
  \tikzexternalenable%
  \tikzsetnextfilename{svd_unfoldings_y}%
  \begin{tikzpicture}[font = \plotfontsize\normalfont]

  \pgfplotstableread[col sep=tab]{graphics/data/svd_unfoldings_Y.dat}\tableY

  \begin{axis}[%
    width            = \textwidth,
    height           = 5.5cm,
    ymode            = log,
    xmin             = 0,
    xlabel           = {Index},
    ylabel           = {$\sigma_i / \sigma_1$},
    grid             = major,
  ]

    \addplot[matlabblue, solid, line width=1.5pt,
             mark=*, mark options={solid}, mark repeat=7, mark size=1.5pt]
      table[x=Index, y=Mode1] {\tableY};

    \addplot[matlaborange, solid, line width=1.5pt,
             mark=square*, mark options={solid}, mark repeat=8, mark size=1.5pt]
      table[x=Index, y=Mode2] {\tableY};

    \addplot[matlabgreen, solid, line width=1.5pt,
             mark=triangle*, mark options={solid}, mark repeat=5, mark size=1.8pt]
      table[x=Index, y=Mode3] {\tableY};

  \end{axis}
\end{tikzpicture}%
  \tikzexternaldisable%

      \end{subfigure}              
      \vspace{0.2cm}

      %
  \tikzexternalenable%
  \tikzsetnextfilename{svd_unfoldings_legend}%
  \begin{tikzpicture}[font = \plotfontsize\normalfont]
  \begin{axis}[%
    hide axis,
    scale only axis,
    width          = .8\linewidth,
    height         = .2\linewidth,
    xmin           = 0,
    xmax           = 1,
    ymin           = 0,
    ymax           = 1,
    legend columns = 3,
    legend style   = {
      at     = {(0,0)},
      anchor = center,
      /tikz/every even column/.append style = {column sep = 0.5cm}}
  ]

    \addplot[matlabblue, solid, line width=1.5pt,
             mark=*, mark options={solid}, mark size=1.5pt] coordinates {(0,0)};
    \addlegendentry{mode-1}

    \addplot[matlaborange, solid, line width=1.5pt,
             mark=square*, mark options={solid}, mark size=1.5pt] coordinates {(0,0)};
    \addlegendentry{mode-2}

    \addplot[matlabgreen, solid, line width=1.5pt,
             mark=triangle*, mark options={solid}, mark size=1.8pt] coordinates {(0,0)};
    \addlegendentry{mode-3}

  \end{axis}
\end{tikzpicture}%
  \tikzexternaldisable%

      \caption{Normalized singular values $\frac{\sigma_i}{\sigma_1}$ of the mode-$k$ unfoldings ($k=1,2,3$) of the ground truth $\bTX_{\mathrm{true}}$  
  (left) and the noisy blurred data $\bTY$ (right). The rapid decay across all modes illustrates the low multilinear rank
   structure exploited by the Tucker-based solvers.}             
      \label{fig:svd_unfolding}    
  \end{figure} 
The effectiveness of the Tucker format hinges on the data being compressible, i.e., admitting accurate low multilinear rank approximations. In an idealized noiseless setting, the data $\opL(\bTX_{\mathrm{true}})$ would be exactly low-rank, meaning that storage and computation scale with the small multilinear ranks $(r_1, r_2, r_3)$ rather than the full tensor dimensions $n_1 \times n_2  \times n_3$ an enormous reduction in both memory and computational cost. In practice, measurements are corrupted by noise, so the observed data $\bTY$ are only \emph{approximately} low-rank. Nevertheless, the underlying low-dimensional structure persists: the singular values of the mode-$k$ unfoldings still decay rapidly, as shown in ~\Cref{fig:svd_unfolding} for both $\bTX_{\mathrm{true}}$ and $\bTY$, confirming that a low-rank Tucker approximation captures the essential information while discarding noise-dominated components. The price of noise is  modest: we must retain slightly larger ranks than in the exact case to obtain an accurate approximation.
We exploit this property by compressing the right-hand  side $\bTY$ via a truncated HOSVD with tolerance $\tau_{\mathrm{outer}}$, yielding a Tucker tensor $\boldsymbol{\mathcal{B}}=\llbracket\boldsymbol{\mathcal{G}};\bU_1,\bU_2,\bU_3\rrbracket$ of multilinear rank $(r_1,r_2,r_3)$ satisfying  $\lVert\boldsymbol{\mathcal{Y}} - \boldsymbol{\mathcal{B}}\rVert_{\rm F} \le \tau_{\mathrm{outer}}                   \lVert\boldsymbol{\mathcal{Y}}\rVert_{\rm F}$. 
This compression introduces an additional perturbation to the right-hand side.
Since $\bTY = \opL(\bTX_{\mathrm{true}}) + \boldsymbol{\eta}$, the compressed tensor $\bTB$ satisfies 
$
\opL(\bTX) \;=\; \bTB,                                     \qquad                                                       \bTB \;=\; \opL(\bTX_{\mathrm{true}})                        + \boldsymbol{\eta}                                               + \boldsymbol{\eta}_{\mathrm{trunc}}$,          
where $\boldsymbol{\eta}$ is the measurement noise and $\boldsymbol{\eta}_{\mathrm{trunc}} = \bTB - \bTY$ is the truncation error introduced by the HOSVD compression, satisfying $\lVert\boldsymbol{\eta}_{\mathrm{trunc}}\rVert_{\fro} \le \tau_{\mathrm{outer}}\lVert\bTY\rVert_{\fro}$. For moderate tolerances $\tau_{\mathrm{outer}}$, this truncation error is small relative to the measurement noise, and because the HOSVD discards components associated with the smallest singular values, the compression preferentially removes high-frequency noise rather than signal. As a result, the low-rank compression has a regularizing effect, suppressing noise before the iterative solver begins. 
To the best of our knowledge, this work is the first to leverage tensor data compression as an implicit regularization mechanism within the iterative solution of inverse problems. To illustrate this effect, we run the MLN-Tucker sGMRES method for 100 iterations with three compression tolerances, $\mathrm{tol} = 0.01, 0.03, 0.05$, and display the RRE versus iteration count in~\Cref{fig:deblurring_res_comp}.  The results show that tighter compression tolerances retain more noise in the right-hand side and lead to earlier onset of semi-convergence, whereas a coarser tolerance ($\mathrm{tol} = 0.05$) effectively suppresses high-frequency noise components and yields a monotonically decreasing RRE, confirming the regularizing role of the truncation.

A key benefit of performing all computations entirely in Tucker format is the substantial reduction in memory requirements.  Storing the full dense tensor
$\bTY\in\R^{128\times 128\times 64}$ requires
$N = 1{,}048{,}576$
floating-point numbers.  By contrast, with $\mathrm{tol} = 0.05$
we obtain a Tucker tensor
$\bTB = \llbracket\bTG;\,\bU_1,\bU_2,\bU_3\rrbracket$ of multilinear
rank $(80,57,3)$, which requires storing only $31{,}408$ floating-point numbers.

Although compression of the right-hand side already provides implicit regularization, the severe ill-conditioning of the blurring operator $\opL$~\cite{Han10} may require additional explicit regularization to fully suppress noise amplification.  To this end, we propose an adaptive projected Tikhonov regularization step directly inside the GMRES iteration, following the hybrid projection approach of~\cite{CalGR99, GazNR15} that we describe below.

\emph{Projected Tikhonov for RHOSVD-Tucker GMRES.}
In the deterministic Tucker GMRES setting
(\Cref{alg:tucker_gmres}), the Arnoldi process produces the relation
$\opL(\bTV_k) = \bTV_{k+1}\overline{\bH}_k$, where
$\overline{\bH}_k\in\R^{(k+1)\times k}$ is upper Hessenberg and $\bTV_1,\dots,\bTV_{k+1}$ are Tucker tensors that are orthonormal in the Frobenius inner product, i.e., $\langle \bTV_i, \bTV_j \rangle_{\fro} = \delta_{ij}$, where $\delta_{ij}$ denotes the Kronecker delta.
In the GMRES algorithm, the coefficient vector $\by_k$ is obtained by minimizing $\lVert\overline{\bH}_k\by -
  \beta\be_1\rVert_2$ (cf.~\eqref{eq:LS:y=Hy-betae1}).  For ill-posed problems this least-squares solution is dominated by noise, so we replace it
  with the Tikhonov-regularized projected problem
\begin{equation}\label{eq:proj_tik}
    \by_k^{(\lambda)}=
    \argmin_{\by\in\R^k}
\bigl\lVert\overline{\bH}_k\by - \beta\be_1\bigr\rVert_2^2 + \lambda_k^2\lVert\by\rVert_2^2.
\end{equation}
Following the hybrid projection framework~\cite{Han10}, regularization is applied to the small $k \times k$ projected system rather than to the original $N \times N$ problem, making the selection of the regularization parameter $\lambda_k$ computationally inexpensive. The approximate tensor solution at iteration $k$
is then assembled as
$ \bTX_k
\;=\;  \texttt{RoundSum}\!\biggl(\sum_{i=1}^k (y_k)_i\,\bTV_i\biggr).$

\emph{Projected Tikhonov for the sketched variants.}
For the sketched methods (RHOSVD-Tucker sGMRES and MLN-Tucker sGMRES),
the Hessenberg matrix is not available because the basis is only
partially orthogonalised. Instead, both algorithms maintain a sketch
matrix $\bM_k\in\R^{s\times k}$ whose $i$-th column is the sketch of
$\opL(\bTV_i)$, i.e.\
$\bM_k = [\bS\,\mathrm{vec}(\opL(\bTV_1)),\dots,
           \bS\,\mathrm{vec}(\opL(\bTV_k))]$,
where $\bS$ is the Khatri--Rao sketching operator.  The regularised
projected problem becomes
\begin{equation}\label{eq:proj_tik_sketched}
  \by_k^{(\lambda)}
  \;=\;
  \argmin_{\by\in\R^k}
  \bigl\lVert\bM_k\by - \bS\,\mathrm{vec}(\bTB)\bigr\rVert_2^2
  + \lambda_k^2\lVert\by\rVert_2^2.
\end{equation}
For MLN-Tucker sGMRES, the sketch
$\bS\,\mathrm{vec}(\opL(\bTV_i))$ is obtained without additional cost
by extracting the diagonal of the Kronecker-sketched core tensor
already computed during the Nystr\"om compression step (see
\Cref{sec: MLN-Tucker sGMRES}).

In both cases, we follow an adaptive approach where the regularization parameter $\lambda_k$ is selected independently at each iteration via generalized cross validation (GCV)~\cite{GolHW79}.~\Cref{fig:deblurring_res_comp}(a) displays the RRE over 100 iterations.~\Cref{fig:deblurring_res_comp}(b) illustrates the effect of compression and regularization. It is clear that in the absence of regularization and compression, the error propagates into the computed solution leading to semiconvergence. 
  
~\Cref{fig:reconstructions}(a)--(c) shows the reconstructions obtained by the unregularised RHOSVD GMRES, RHOSVD sGMRES, and MLN sGMRES methods, while~\Cref{fig:reconstructions}(d)--(f) displays the corresponding           
  reconstructions when right-hand side compression and GCV regularization are
  applied.  
\begin{figure}[htbp]
    \centering
    \begin{subfigure}[b]{0.48\textwidth}
        \centering
  \tikzexternalenable%
  \tikzsetnextfilename{deblurring_residuals_unreg}%
  \begin{tikzpicture}

    \pgfplotstableread[col sep=tab]{graphics/data/deblurring_unreg_MLNTuckersGMRES_history_k2.dat}\tableMLNsGMRES
    \pgfplotstableread[col sep=tab]{graphics/data/deblurring_unreg_MLNTuckersGMRESMemoryEfficient_history_k2.dat}\tableMLNsGMRESMemoryEfficient
    \pgfplotstableread[col sep=tab]{graphics/data/deblurring_unreg_RHOSVDTuckersGMRES_history_k2.dat}\tablesGMRES
    \pgfplotstableread[col sep=tab]{graphics/data/deblurring_unreg_RHOSVDTuckerGMRES_history_k2.dat}\tableGMRES

    \begin{axis}[
        width             = \textwidth,
        height            = 5.5cm,
        ymode             = log,
        ytick             = {0.2, 0.4, 0.6, 0.8, 1.0},
        xmin              = 0,
        xmax              = 100,
        ymin              = 0.2,
        ymax              = 1.0,
        minor y tick num  = 5,
        log ticks with fixed point,
        xlabel            = {Iteration},
        ylabel            = {RRE},
        grid              = major,
    ]
    
    \addplot[MLNsGMRES, mark repeat=7]
        table[x=Iteration, y=RRE] {\tableMLNsGMRES};

    \addplot[MLNsGMRESMemoryEfficient, mark repeat=6]
        table[x=Iteration, y=RRE] {\tableMLNsGMRESMemoryEfficient};

    \addplot[sGMRES, mark repeat=5]
        table[x=Iteration, y=RRE] {\tablesGMRES};

    \addplot[GMRES, mark repeat=8]
        table[x=Iteration, y=RRE] {\tableGMRES};
    
    \end{axis}
\end{tikzpicture}%
  \tikzexternaldisable%

        \caption{No compression, no regularization}
        \label{fig:deblurring_res_unreg}
    \end{subfigure}
    \hfill
    \begin{subfigure}[b]{0.48\textwidth}
        \centering
  \tikzexternalenable%
  \tikzsetnextfilename{deblurring_residuals_reg}%
  \begin{tikzpicture}

    \pgfplotstableread[col sep=tab]{graphics/data/deblurring_regcomp_MLNTuckersGMRES_history_k2.dat}\tableMLNsGMRES
    \pgfplotstableread[col sep=tab]{graphics/data/deblurring_regcomp_MLNTuckersGMRESMemoryEfficient_history_k2.dat}\tableMLNsGMRESMemoryEfficient
    \pgfplotstableread[col sep=tab]{graphics/data/deblurring_regcomp_RHOSVDTuckersGMRES_history_k2.dat}\tablesGMRES
    \pgfplotstableread[col sep=tab]{graphics/data/deblurring_regcomp_RHOSVDTuckerGMRES_history_k2.dat}\tableGMRES

    \begin{axis}[
        width             = \textwidth,
        height            = 5.5cm,
        ymode             = log,
        ytick             ={0.3, 0.4, 0.5},
        xmin              = 0,
        xmax              = 100,
        ymin              = 0.25,
        ymax              = 0.5,
        minor y tick num  = 5,
        log ticks with fixed point,
        xlabel            = {Iteration},
        ylabel            = {RRE},
        grid              = major,
    ]
    
    \addplot[MLNsGMRES, mark repeat=7]
        table[x=Iteration, y=RRE] {\tableMLNsGMRES};

    \addplot[MLNsGMRESMemoryEfficient, mark repeat=6]
        table[x=Iteration, y=RRE] {\tableMLNsGMRESMemoryEfficient};

    \addplot[sGMRES, mark repeat=5]
        table[x=Iteration, y=RRE] {\tablesGMRES};

    \addplot[GMRES, mark repeat=8]
        table[x=Iteration, y=RRE] {\tableGMRES};
    
    \end{axis}
\end{tikzpicture}%
  \tikzexternaldisable%

        \caption{Compression and regularization}
        \label{fig:deblurring_res_reg}
    \end{subfigure}
    
    \vspace{0.2cm}
    
    %
  \tikzexternalenable%
  \tikzsetnextfilename{deblurring_residuals_legend}%
  \begin{tikzpicture}[font = \plotfontsize\normalfont]
  \begin{axis}[%
    hide axis,
    scale only axis,
    width          = .8\linewidth,
    height         = .2\linewidth,
    xmin           = 0,
    xmax           = 1,
    ymin           = 0,
    ymax           = 1,
    legend columns = 2,
    legend style   = {
      at     = {(0,0)},
      anchor = center,
      /tikz/every even column/.append style = {column sep = 0.5cm}}
  ]
    
    \addplot[MLNsGMRES, mark repeat = 7] coordinates {(0,0)};
    \addlegendentry{MLN sGMRES}

    \addplot[MLNsGMRESMemoryEfficient, mark repeat = 6] coordinates {(0,0)};
    \addlegendentry{MLN sGMRES (Memory Efficient)}

    \addplot[sGMRES, mark repeat = 5] coordinates {(0,0)};
    \addlegendentry{RHOSVD sGMRES}

    \addplot[GMRES, mark repeat = 8] coordinates {(0,0)};
    \addlegendentry{RHOSVD GMRES}
    
  \end{axis}
\end{tikzpicture}%
  \tikzexternaldisable%

    \caption{Convergence histories comparing RRE across the unregularised and regularised setups. All variants truncated at $k_{\text{trunc}} = 2 $.}
    \label{fig:deblurring_residuals_methods}
\end{figure}
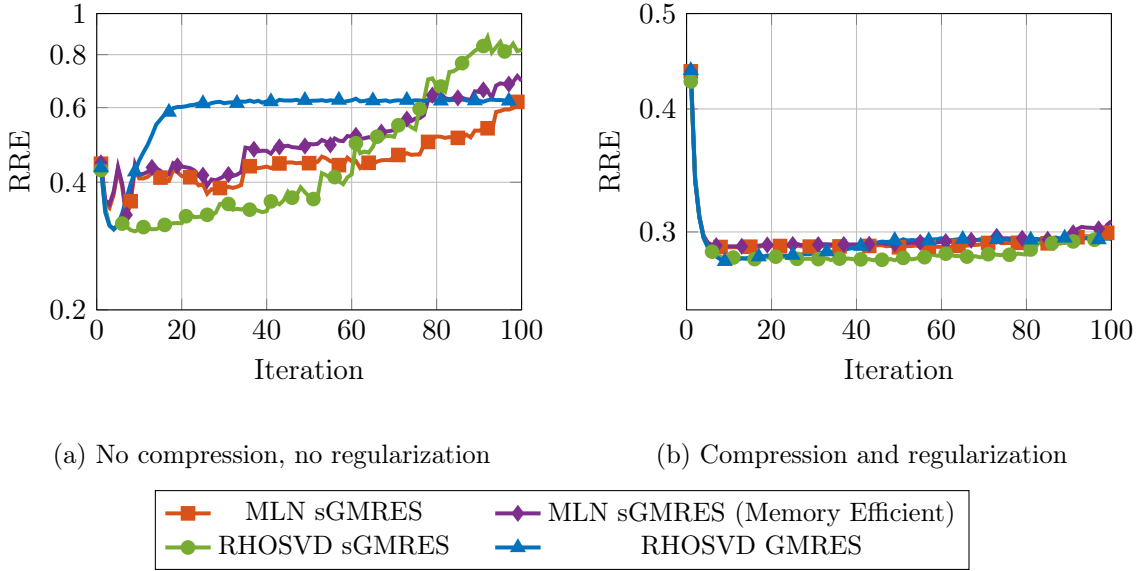
\begin{figure}[t]
\centering

\begin{subfigure}[t]{0.24\textwidth}
\centering
\includegraphics[width=\textwidth]{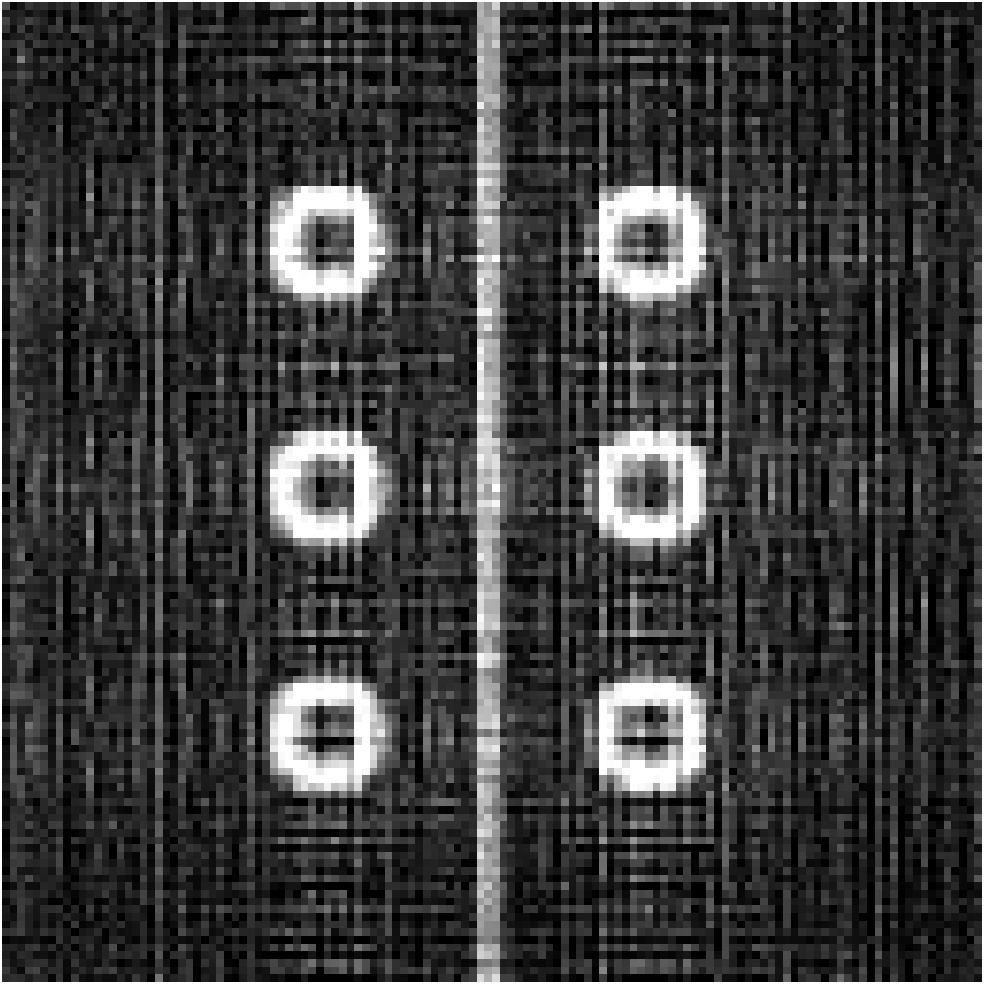}
\caption{}
\end{subfigure}
\hfill
\begin{subfigure}[t]{0.24\textwidth}
\centering
\includegraphics[width=\textwidth]{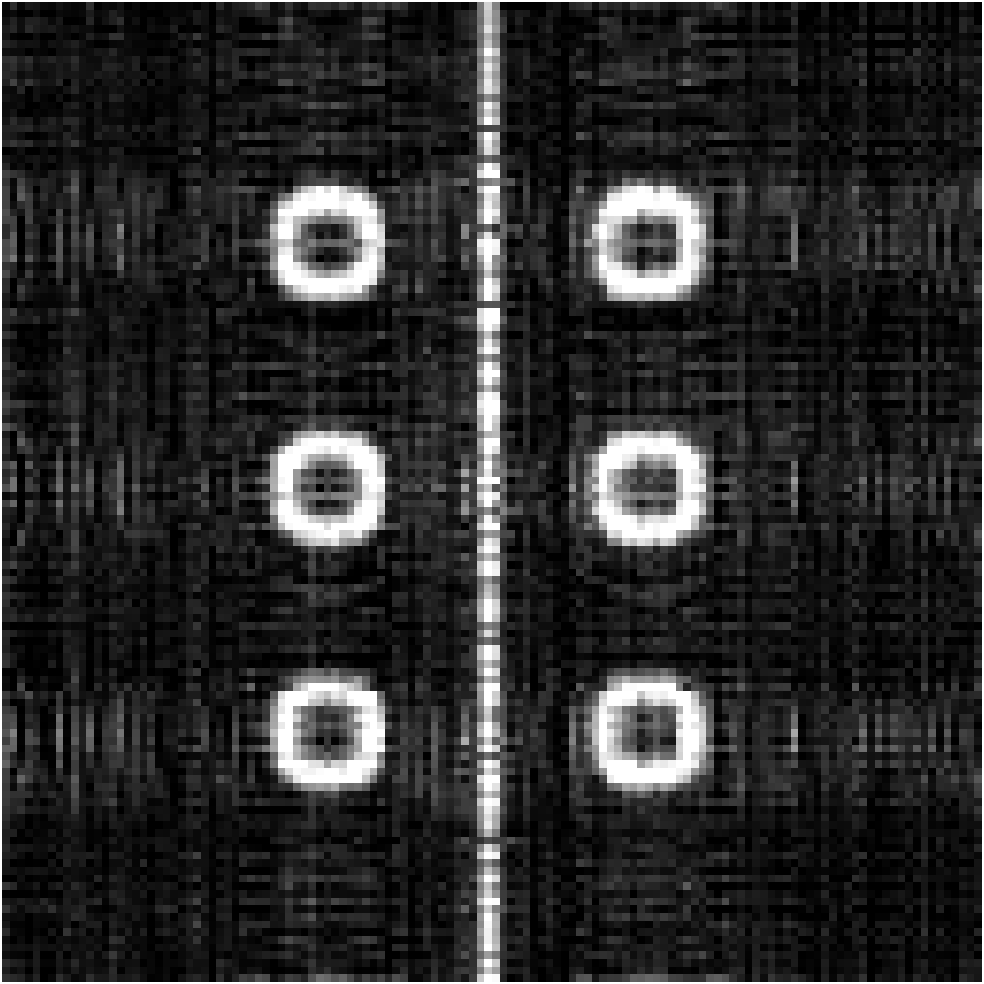}
\caption{}
\end{subfigure}
\hfill
\begin{subfigure}[t]{0.24\textwidth}
\centering
\includegraphics[width=\textwidth]{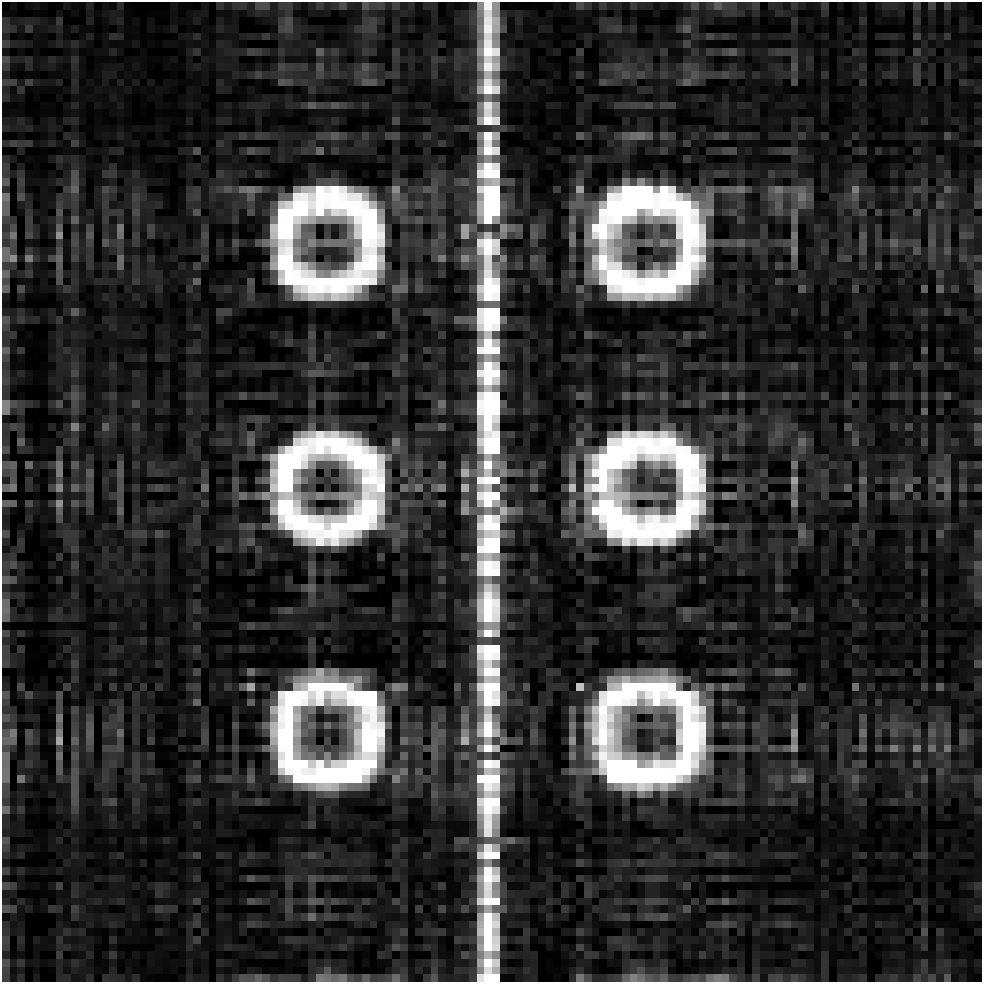}
\caption{}
\end{subfigure}
\hfill
\begin{subfigure}[t]{0.24\textwidth}
\centering
\includegraphics[width=\textwidth]{graphics/data/deblurring_unreg_memory_efficient_RHOSVDTuckersGMRES.png}
\caption{}
\end{subfigure}

\vspace{0.2cm}

\begin{subfigure}[t]{0.24\textwidth}
\centering
\includegraphics[width=\textwidth]{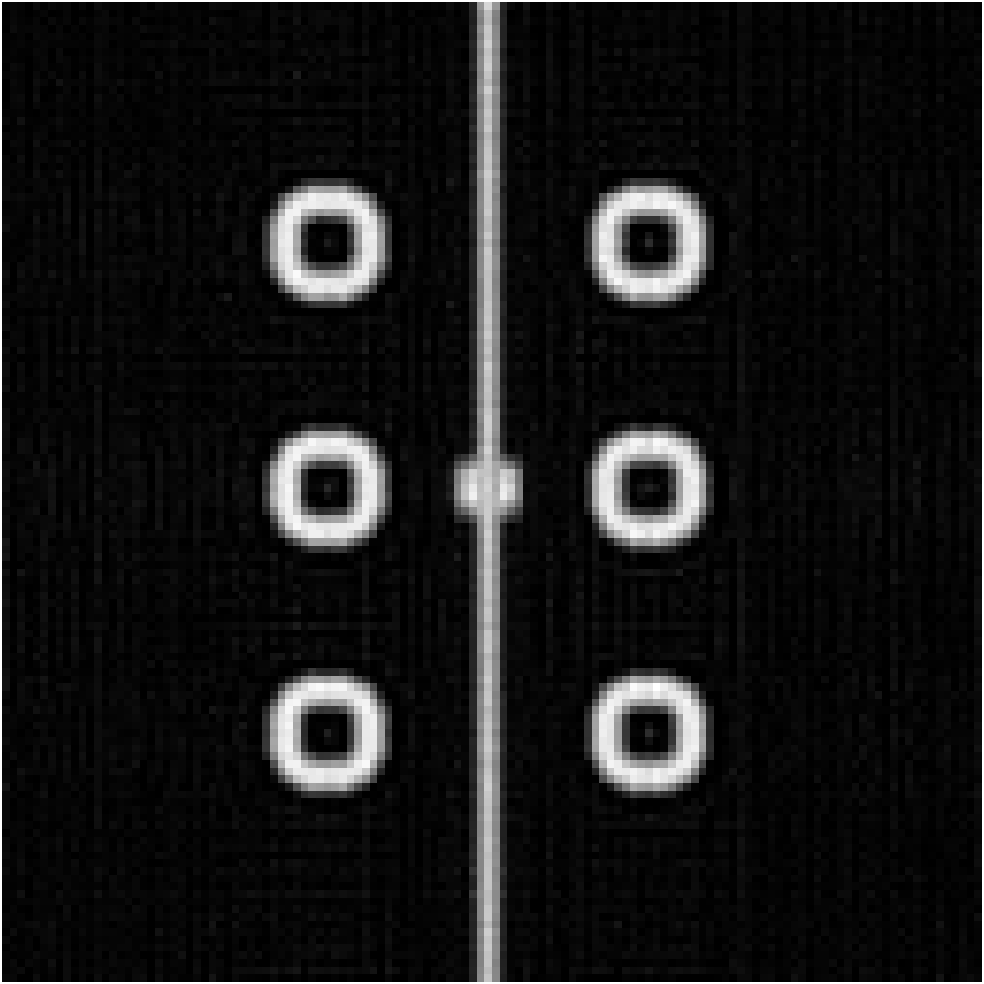}
\caption{}
\end{subfigure}
\hfill
\begin{subfigure}[t]{0.24\textwidth}
\centering
\includegraphics[width=\textwidth]{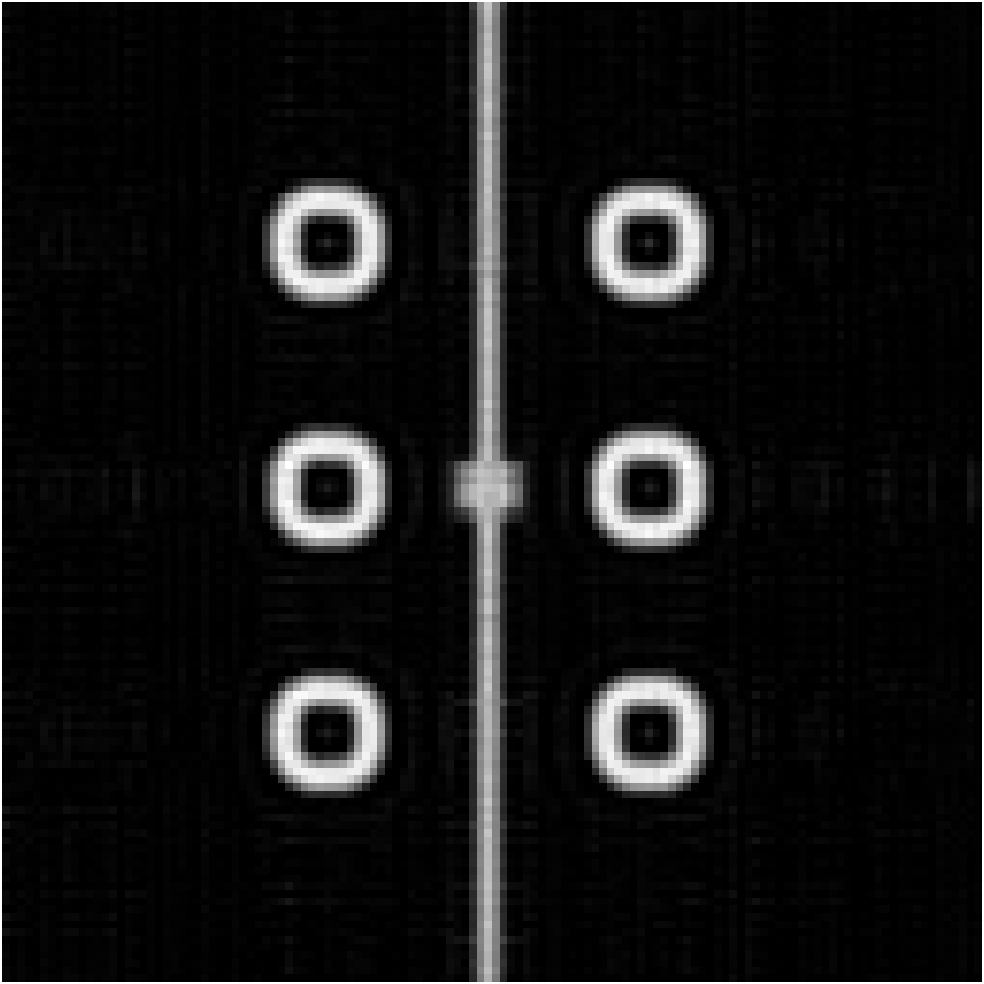}
\caption{}
\end{subfigure}
\hfill
\begin{subfigure}[t]{0.24\textwidth}
\centering
\includegraphics[width=\textwidth]{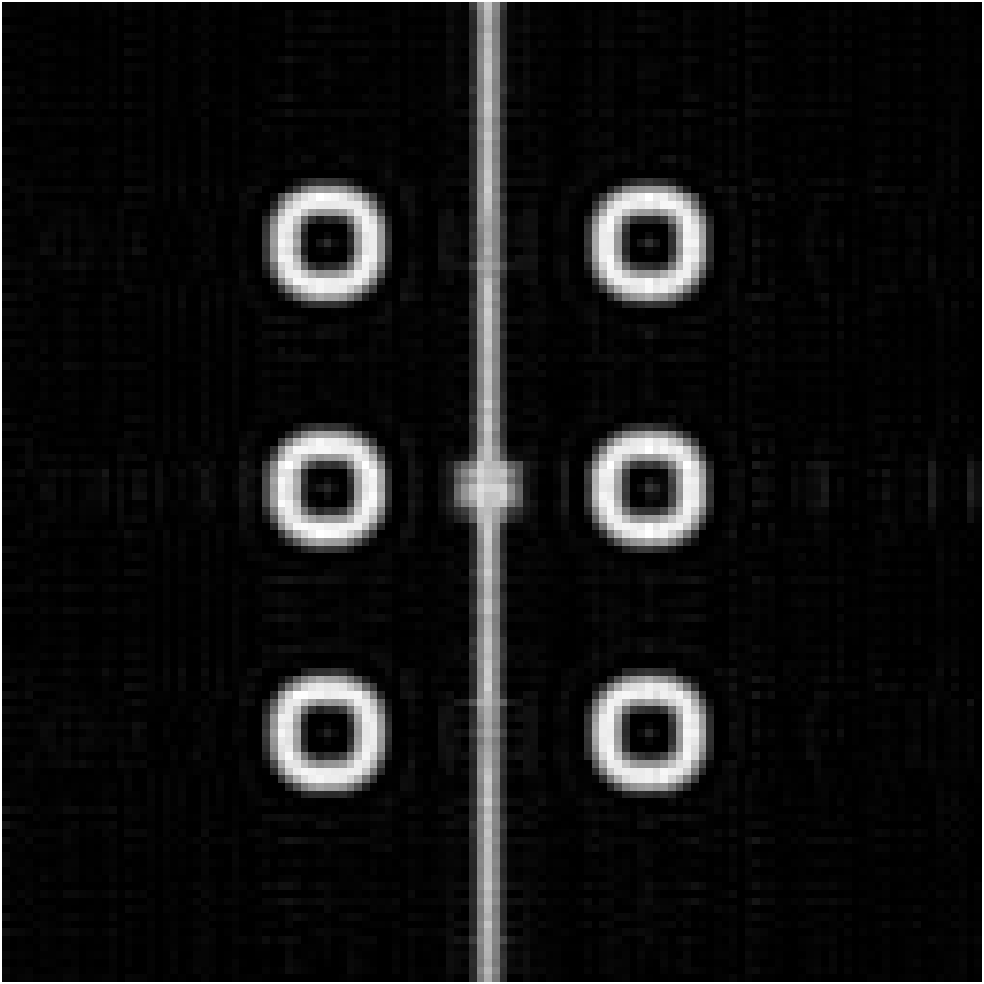}
\caption{}
\end{subfigure}
\hfill
\begin{subfigure}[t]{0.24\textwidth}
\centering
\includegraphics[width=\textwidth]{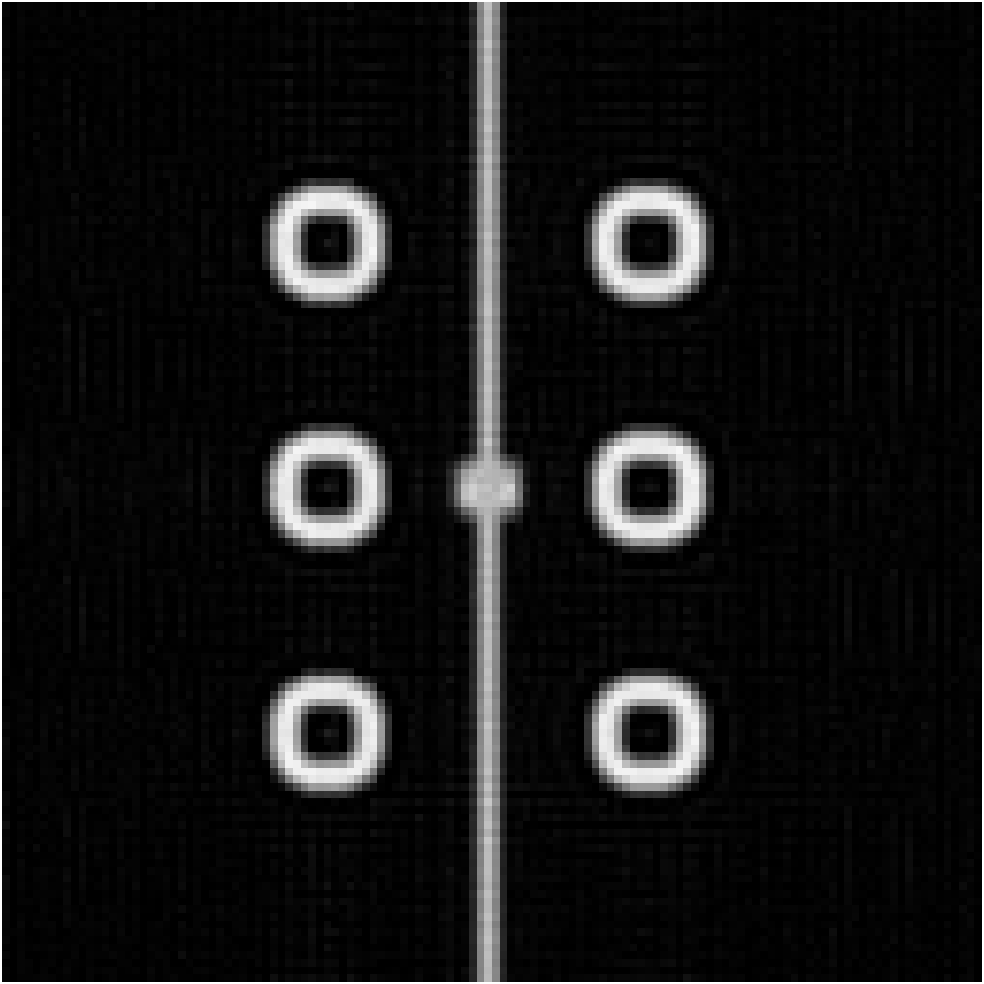}
\caption{}
\end{subfigure}

\caption{Image deblurring reconstructions.
Top row: unregularised (a)~RHOSVD sGMRES, (b)~MLN sGMRES, (c)~MLN sGMRES (\texttt{save\_memory}), (d)~RHOSVD GMRES.
Bottom row: compressed RHS + GCV regularization (e)~RHOSVD sGMRES, (f)~MLN sGMRES, (g)~MLN sGMRES (\texttt{save\_memory}), (h)~RHOSVD GMRES.}
\label{fig:reconstructions}
\end{figure}

\section{Conclusions}
\label{sec:conclusions}
We have presented two novel randomized iterative solvers, namely, the RHOSVD-Tucker sGMRES and the MLN-Tucker sGMRES for large-scale tensor equations expressed as a sum of separable operators. Both methods operate entirely in the Tucker format, and work by combining sketched GMRES and randomized compression.

RHOSVD-Tucker sGMRES extends Khatri–Rao sketching with adaptive rank estimation and replaces full Arnoldi orthogonalization with a short-recurrence partial orthogonalization that controls the multilinear ranks of the basis tensors. In contrast, MLN-Tucker sGMRES leverages a single-pass and streamable multilinear Nystr\"om approximation: by fixing the rank a priori and discarding basis tensors after sketching, it avoids storing the full Krylov basis and, under low-rank assumption in the solution, allows to operate in compressed sketch spaces, leading to substantial savings in memory. Moreover, the cost is reduced as the Khatri–Rao sketch required by sGMRES is obtained as a by-product of the Nystr\"om compression.

Numerical experiments on representative 3D problems, including Poisson and convection-diffusion equations, as well as an image deblurring inverse problem, show that the proposed methods achieve accurate solutions while reducing memory usage and computational effort by orders of magnitude compared to dense approaches. For the deblurring problem, compressing the data into Tucker format serves a dual purpose: it makes an otherwise intractable large-scale problem computationally accessible, and, combined with an adaptive projected Tikhonov penalty selected via GCV, acts as an implicit regularizer that yields stable reconstructions without prior knowledge of the noise level. To our knowledge, this is the first work to exploit compressed tensor decompositions as both a regularization mechanism and an enabler of tractable computation within iterative Krylov solvers for inverse problems.
\section*{Acknowledgments}
MP acknowledges support from the National Science Foundation (NSF) under Grant No. DMS-2410699. Any opinions, findings, conclusions, or recommendations expressed in this material are those of the authors and do not necessarily reflect the views of the National Science Foundation. RS and MP acknowledge Advanced Research Computing at Virginia Tech for providing computational resources and technical support that have contributed to the results reported within this paper. URL: https://arc.vt.edu/. RS, MP, and AB acknowledge support from the NSF under Grant No. DMS-1929284 while the authors were in residence at the Institute for Computational and Experimental Research in Mathematics in Providence, RI, during the \emph{Stochastic and Randomized Algorithms in Scientific Computing: Foundations and Applications} semester program, where preliminary discussions of this work took place. 

\addcontentsline{toc}{section}{References}
\bibliographystyle{plainurl}
\bibliography{bibtex/myref}

\end{document}